\documentclass[a4paper]{article}
\usepackage{amsmath,amsthm,amssymb,bm,1-dimensional_sum}
\usepackage{mathscinet}
\usepackage{amscd}
\usepackage[all]{xy}

\begin{document}
\title{Demazure crystals and tensor products of perfect Kirillov-Reshetikhin crystals with various levels}
\author{Katsuyuki Naoi}
\date{}

\maketitle

\begin{abstract}
  In this paper, we study a tensor product of perfect Kirillov-Reshetikhin crystals (KR crystals for short)
  whose levels are not necessarily equal.
  We show that, by tensoring with a certain highest weight element,
  such a crystal becomes isomorphic as a full subgraph
  to a certain disjoint union of Demazure crystals contained in a tensor product of highest weight crystals.
  Moreover, we show that this isomorphism preserves their $\Z$-gradings, 
  where the $\Z$-grading on the tensor product of KR crystals is given by the energy function, 
  and that on the other side is given by the minus of the action of the degree operator.
\end{abstract}

\section{Introduction}

Crystal bases $B(\gL)$ introduced by Kashiwara \cite{MR1115118}
can be viewed as bases at $q=0$ of highest weight modules $V(\gL)$ of the quantized enveloping algebra $U_q(\fg)$
associated with a Kac-Moody Lie algebra $\fg$.
Crystal bases reflect the internal structures of the modules, and are powerful combinatorial tools for studying them.

Crystal bases are also useful for studying certain subspaces of $V(\gL)$.
For a Weyl group element $w$, the Demazure module $V_w(\gL)$, which is a module of a Borel subalgebra,
is defined 
by the submodule of $V(\gL)$ generated by the extremal weight space $V(\gL)_{w\gL}$.
Kashiwara showed in \cite{MR1240605} that there exists a certain subset $B_w(\gL) \subseteq B(\gL)$ 
which is, in a suitable sense, a crystal basis of $V_w(\gL)$.
The subset $B_w(\gL)$ is called the Demazure crystal.
Using Demazure crystals, he gave a new proof of the character formula for Demazure modules in the article,
which expresses the character using the Demazure operators (see \cite{MR1240605} or 
Section \ref{section:Demazure_crystals} of the present article).

When $\fg$ is an affine Kac-Moody Lie algebra, there is another class of modules having crystal bases 
called Kirillov-Reshetikhin modules $W^{r,\ell}$ (KR modules for short), where $r$ is a node in the classical Dynkin diagram
and $\ell$ is a positive integer.
KR modules are finite-dimensional irreducible $U_q'(\fg)$-modules, 
where $U_q'(\fg)$ is the quantum affine
algebra without the degree operator. 
At least when $\fg$ is nonexceptional, it was proved that 
every $W^{r,\ell}$ has a crystal basis $B^{r,\ell}$ \cite{MR1194953,MR2389790,MR2403558},
which is called the Kirillov-Reshetikhin (KR) crystal.

Demazure crystals and KR crystals are known to have strong relations, 
and the study of the relationship between them has been the subject of many articles.
For example, see \cite{MR2301245,MR2115972,MR1643999,MR1620507,MR2855081,ST}.

Among these articles, \cite{ST} by Schilling and Tingley is quite important for the present article.
In the article, they studied
a tensor product of perfect KR crystals for nonexceptional $\fg$ whose levels are all the same
(perfectness is a technical condition for a finite $U_q'(\fg)$-crystal
which allows one to use the crystal to construct highest weight crystals, see \cite{MR1187560} or 
Subsection \ref{subsection:Perfect_KR_crystals} of the present article).
They proved that, by tensoring with a certain highest weight element, 
such a crystal becomes isomorphic to a certain Demazure crystal as a full subgraph.
Moreover, they also showed that this isomorphism preserves their $\Z$-gradings.
Here, the tensor product of perfect KR crystals is $\Z$-graded by the energy function, which is a certain $\Z$-function
defined in a combinatorial way, and the $\Z$-grading of the Demazure crystal is given by 
the minus of the action of the degree operator.
Since the Demazure crystal has a character formula as stated above,
these results imply that the weight sum with the energy function 
of the tensor product of perfect KR crystals (with same levels) can be expressed by the Demazure character formula.

The aim of this article is to generalize the above results to a tensor product of perfect KR 
crystals whose levels are not necessarily equal.
In this case the tensor product of perfect KR crystals, tensored with a highest weight element, 
is no longer isomorphic to a single Demazure crystal.
We show in this article, however, that
it is isomorphic to a certain disjoint union of Demazure crystals contained in a tensor product of $B(\gL)$'s,
and that this isomorphism also preserves their $\Z$-gradings.

Before stating our results, we prepare some notation.
For a crystal $B$ and a Dynkin automorphism $\tau$, we define a new crystal 
$\tilde{\tau}(B) = \{ \tilde{\tau}(b)\mid b \in B\}$ whose weight function is $\wt\big(\tilde{\tau}(b)\big) =
\tau\big(\wt(b)\big)$
and Kashiwara operators are 
\[ e_i \tilde{\tau}(b) = \tilde{\tau}\big(e_{\tau^{-1}(i)} b\big), \ \ \ f_i \tilde{\tau}(b) = \tilde{\tau}
   \big(f_{\tau^{-1}(i)}b\big).
\]
Let $S$ be a subset of $B$, and $w$ a Weyl group element with a reduced expression $w =s_{i_k} \cdots s_{i_1}$.
We denote by $\cF_{w\tau}(S)$ the subset of $\tilde{\tau}(B)$ defined by
\[ \cF_{w\tau}(S) = \bigcup_{j_1,\dots,j_k \ge 0} f_{i_k}^{j_k} \cdots f_{i_1}^{j_1} \tilde{\tau}(S) \setminus \{0\}.
\]
For a dominant integral weight $\gL$, we denote by $u_{\gL}$ the highest weight element of $B(\gL)$.

Now let us mention our results.
Assume that $\fg$ is of nonexceptional type,
and let $B^{r_1,c_{r_1}\ell_1}, \dots, B^{r_p,c_{r_p}\ell_p}$ be perfect KR crystals.
Here, $c_r$ is a particular constant which ensures the perfectness for the KR crystal $B^{r,c_r\ell}$.
We assume $\ell_1 \le \ell_2 \le \cdots \le \ell_p$, and
put $\mu_j = c_{r_j}w_0(\varpi_{r_j})$ for $1 \le j \le p$, where $w_0$ is the longest element of the Weyl group of
the simple Lie subalgebra $\fg_0 \subseteq \fg$ corresponding to the classical Dynkin diagram, and 
$\varpi_{r}$ are the fundamental weights of $\fg_0$.
Then the following theorem is proved, which is the main theorem of the present article (Theorem \ref{Thm:Main_Theorem}):

\begin{Thm}\label{Thm:intro}
  There exists an isomorphism
  \begin{align*} 
    u_{\ell_p\gL_0} \otimes &B^{r_p,c_{r_p}\ell_p} \otimes \cdots \otimes B^{r_1,c_{r_1}\ell_1} \\ \stackrel{\sim}{\to}
    &\cF_{t_{\mu_p}}\Big( u_{(\ell_p - \ell_{p-1})\gL_0} \otimes \dots \otimes \cF_{t_{\mu_2}}\Big(u_{(\ell_2 - \ell_1)\gL_0} 
    \otimes \cF_{t_{\mu_1}}(u_{\ell_1\gL_0})\Big)\!\cdots\!\Big)
  \end{align*}
  of full subgraphs, where $t_{\mu}$ denotes the translation and $\gL_0$ denotes the fundamental weight of $\fg$. 
  Moreover, this isomorphism preserves the $\Z$-gradings up to a shift.
\end{Thm}

Using the combinatorial excellent filtration theorem \cite{MR1887117,MR1987017},
it is easy to see that the right hand side of the above isomorphism is a disjoint union of Demazure crystals.
Then similarly as a Demazure crystal, the weight sum of the right hand side
can be expressed using Demazure operators.
Hence, we obtain the following corollary (Corollary \ref{Cor:Main_cor}),
where we set $B = B^{r_p,c_{r_p}\ell_p} \otimes \cdots \otimes B^{r_1,c_{r_1}\ell_1}$:

\begin{Cor}\label{Cor:introA}
  Let $\aff:P_\cl \to P$ denote the canonical section of the projection from affine weight lattice $P$ 
  to the classical weight lattice $P_{\cl}$. Then we have
  \begin{align*}
    e^{\ell_p\gL_0+ C_B\gd} &\sum_{b \in B} e^{\aff\circ\wt(b) - \gd D(b)}\\ 
    &= D_{t_{\mu_p}}\left( e^{(\ell_p-\ell_{p-1})\gL_0} \cdots D_{t_{\mu_2}}
    \left( e^{(\ell_2 -\ell_1)\gL_0} \cdot D_{t_{\mu_1}}(e^{\ell_1\gL_0})\right)\!\cdots\! \right)
  \end{align*}
  for some constant $C_B$,
  where $D_{t_\mu}$ is the Demazure operator associated with $t_\mu$ 
  {\normalfont(}see Section \ref{section:Demazure_crystals}{\normalfont)}, and $D\colon B \to \Z$ is the energy function.  
\end{Cor}

Let $X(B, \mu ,q)$ denote the one-dimensional sum \cite{MR1903978,MR1745263} associated with
the crystal $B$ and a dominant integral weight $\mu$ of $\fg_0$.
Then the above corollary is equivalent to the following (Corollary \ref{Cor:X}):

\begin{Cor} \label{Cor:introB}
  Let $P_0^+$ be the set of dominant integral weights of $\fg_0$ 
  and $\ch V_{\fg_0}(\mu)$ the character of the irreducible $\fg_0$-module with highest weight $\mu$.
  Then we have
  \begin{align}\label{eq:Cor_intro}
    q^{-C_B}\sum_{\mu \in P_0^+}& X(B, \mu,q) \ch V_{\fg_0}(\mu) \\ 
    &= e^{-\ell_p\gL_0}D_{t_{\mu_p}}\left( e^{(\ell_p-\ell_{p-1})\gL_0} \cdots D_{t_{\mu_2}}
    \left( e^{(\ell_2 -\ell_1)\gL_0} \cdot D_{t_{\mu_1}}(e^{\ell_1\gL_0})\right)\!\cdots\! \right), \nonumber
  \end{align}
  where we set $q = e^{-\gd}$.
\end{Cor}

Corollary \ref{Cor:introB} has an important application (and in fact this is one of the main motivations of this work).
The $X=M$ conjecture presented in \cite{MR1903978,MR1745263} asserts that a one-dimensional
sum is equal to a fermionic form which is defined as a generating function of some combinatorial objects called 
rigged configurations.
In \cite{MR2964614}, it is proved that when $\fg$ is of type $A_n^{(1)}, D_n^{(1)}$ or $E_n^{(1)}$, 
the fermionic forms also satisfy a similar equation as (\ref{eq:Cor_intro}).
Then we can prove the $X=M$ conjecture in the cases $\fg = A_n^{(1)}, D_n^{(1)}$ from these equations
(for details see \cite{MR2964614}).

The plan of this article is as follows.
In Section 2, we fix basic notation used in the article.
In Section 3, we briefly review the definition of crystals, and
in Section 4, we review the results on Demazure crystals.
In Section 5, we review the results on KR crystals, and construct the isomorphism in Theorem \ref{Thm:intro}.
In Section 6, we review the definition and some results on the energy functions,
and finally in Section 7, we show that the isomorphism constructed in Section 5 preserves the $\Z$-gradings,
which completes the proof of Theorem \ref{Thm:intro}.

\ \\
\textbf{Acknowledgements:} 
The author is very grateful to M.\ Okado for answering many questions and providing a lot of references.
Without his help, this paper could not have been written.
He also thank R.\ Kodera for reading the manuscript very carefully and pointing out a lot of errors, S.\ Naito and Y.\ Saito
for some helpful comments, and A.\ Schilling for sending him her preprint.
This work was supported by World Premier International Research Center Initiative (WPI Initiative), MEXT, Japan.

\section{Notation and basics}

\subsection{Affine Kac-Moody Lie algebra}\label{subsection:Affine}

Let $\fg$ be a complex affine Kac-Moody Lie algebra with Cartan subalgebra $\fh$, Dynkin node set $I = \{ 0,\dots,n\}$,
Dynkin diagram $\gG$ and Cartan matrix $A = ( a_{ij})_{i,j \in I}$.
In this article, we use the Kac's labeling of nodes of Dynkin diagrams in \cite[Section 4.8]{MR1104219}.
Let $\ga_i\in \fh^{*}$ and $\ga^{\vee}_i \in \fh$ ($i \in I$) be 
the simple roots and the simple coroots respectively, and $\gD \subseteq \fh^*$ the root system of $\fg$.
Let $(a_0,\dots,a_n)$ (resp. $(a_0^{\vee},\dots,a_n^{\vee})$) be the unique sequence of relatively prime positive integers 
satisfying 
\[ \sum_{j \in I} a_{ij}a_j = 0 \ \text{for all $i \in I$} \ (\text{resp.} 
   \sum_{i \in I} a_i^{\vee} a_{ij} = 0 \ \text{for all $j \in I$}).
\]
Let $d \in \fh$ be the degree operator, which is any element satisfying $\langle \ga_i, d \rangle = \gd_{0i}$ for $i \in I$,
$K = \sum_{i \in I} a_i^{\vee}\ga_i^{\vee} \in \fh$ the canonical central element,
$\gd = \sum_{i \in I} a_i \ga_i \in \fh^*$ the null root, and 
$W$ the Weyl group of $\fg$ with simple reflections $\{s_i \mid i \in I\}$.
Let $\ell\colon W \to \Z_{\ge 0}$ be the length function.
Let $Q$ be the root lattice of $\fg$, and $Q^+ = \sum_{i \in I} \Z_{\ge 0} \ga_i$.
In this article we fix a positive integer $N$, and define the weight lattice of $\fg$ by
\begin{equation} \label{eq:def_of_P}
  P = \{ \gl \in \fh^* \mid \langle \gl, \ga_i^{\vee} \rangle \in \Z \ (i \in I), \langle \gl, d \rangle \in N^{-1} \Z \}.
\end{equation}
(In the next subsection, we impose some condition on $N$ so that $P$ is preserved
 by the action of the extended affine Weyl group $\wti{W}$.)
Put $P^+ = \{ \gl \in P \mid \langle \gl, \ga_i^{\vee} \rangle \in \Z_{\ge 0} \ (i \in I)\}$,
and let $\gL_i \in P^+ \ (i \in I)$ be any element satisfying 
\[ \langle \gL_i, \ga_j^{\vee} \rangle = \gd_{ij} \ \ \ \text{for} \ j \in I.
\]
(Note that we do not assume $\langle \gL_i, d \rangle = 0$.)
Then we have $P^+ = \sum_{i \in I} \Z_{\ge 0} \gL_i + N^{-1}\Z\gd$.
For $\gl \in P$, we call the integer $\langle \gl, K \rangle$ the \textit{level} of $\gl$,
and for $\ell \in \Z$ we denote 
$P^\ell = \{\gl \in P \mid \langle \gl, K \rangle = \ell \}$.
Let $( \ , \ )$ be a $W$-invariant symmetric bilinear form on $\fh^*$ satisfying
\[ (\ga_i, \ga_j) = a_i^{\vee}a_i^{-1}a_{ij} \ \text{for} \ i,j\in I, \ \ \ (\ga_i, \gL_0) =\gd_{0i}a_0^{-1} \ \text{for} \ 
   i \in I.
\]

Let $\cl \colon \fh^* \to \fh^*/ \C \gd$ be the canonical projection, and put $P_{\cl} = \cl(P)$, $P_\cl^+ = \cl(P^+)$,
$P_{\cl}^\ell = \cl(P^{\ell})$ for $\ell \in \Z$ and $(P_{\cl}^+)^\ell = P_{\cl}^+ \cap P_{\cl}^{\ell}$.
Since $W$ fixes $\gd$, $W$ acts on $\fh^*/ \C\gd$ and $P_{\cl}$.
For $i \in I$, define $\varpi_i \in P_\cl^0$ by 
$\varpi_i = \cl(\gL_i) - a_i^{\vee}\cl(\gL_0)$.
Note that $\varpi_0 = 0$ and $\varpi_i$ for $i \in I \setminus \{0\}$ satisfies 
\[ \langle \varpi_i, \ga_j^{\vee}\rangle = \gd_{ij} \ \text{for} \ j \in I\setminus \{ 0 \}, \ \ \
   \langle \varpi_i, \ga_0^{\vee}\rangle = -a_i^{\vee}.
\] 
We define $\aff\colon \fh^*/\C \gd \to \fh^*$ by the unique section of $\cl$ satisfying $\langle \aff(\gl), d \rangle = 0$ 
for all $\gl \in \fh^*/ \C \gd$. 
When no confusion is possible, we omit the notation $\cl( * )$ for simplicity.
In particular, we often write $\cl(\gL_i)$ and $\cl(\ga_i)$ simply as $\gL_i$ and $\ga_i$.

Let $I_0 = I \setminus \{ 0 \}$ and $\fg_0 \subseteq \fg$ the simple Lie subalgebra whose Dynkin node set is $I_0$
with Cartan subalgebra $\fh_0 \subseteq \fh$ and Weyl group $W_0 \subseteq W$. 
Let $\varpi_j^{\vee} \in \fh_0$ $(j \in I_0)$ be the unique element satisfying 
\[ \langle \ga_i, \varpi_j^{\vee} \rangle = \gd_{ij} \ \text{for} \ i \in I_0,
\]
which also satisfies 
\begin{equation} \label{eq:pairing}
  \langle \ga_0, \varpi_j^{\vee} \rangle = -a_j/a_0.
\end{equation}
For the notational convenience, we put $\varpi_0^{\vee} = 0$.
Denote by $w_0$ the longest element of $W_0$.
Let $P_0$ denote the weight lattice of $\fg_0$, $P_0^+ \subseteq P_0$ the set of dominant integral weights,
$Q_0 \subseteq P_0$ the root lattice, and $Q_0^+ = \sum_{i \in I_0} \Z_{\ge 0} \ga_i$.
We often identify $P_\cl^0$ and $P_0$ in an obvious way.

The bilinear form $( \ , \ )$ induces a bilinear form on $P_{\cl}^0$, which is also denoted by $( \ , \ )$.
Then we have 
\begin{equation}\label{eq:final}
  (\gl, \varpi_i) = a_i^{\vee}a_i^{-1} \langle \gl, \varpi_i^{\vee} \rangle
\end{equation}
for $i \in I$ and $\gl \in P_\cl^0$.

\subsection{Dynkin automorphisms and extended affine Weyl group}\label{subsection:Dynkin}

As \cite[(6.5.2)]{MR1104219}, we define for $\gl \in P_{\cl}^0$ an endomorphism $t_{\gl}$ of $\fh^*$ by
\begin{equation} \label{eq:translation}
  t_{\gl}(\mu) = \mu + \langle \mu, K \rangle \aff(\gl) - \Big(\big(\mu, \aff(\gl)\big) + 
  \frac{1}{2}\big(\aff(\gl), \aff(\gl)\big)\langle  \mu, K \rangle\Big)\gd.
\end{equation}
The map $\gl \mapsto t_{\gl}$ defines an injective group homomorphism from $P_\cl^0$ to the group of linear 
automorphisms of $\fh^*$ orthogonal with respect to $( \ , \ )$.
Let $c_i = \max\{1,a_i/a_i^{\vee}\}$ for $i \in I_0$,
and define sublattices $M$ and $\wti{M}$ of $P_{\cl}^0$ by
\[ M = \sum_{w \in W_0} \Z w (\ga_0/a_0), \ \ \ \widetilde{M} = \bigoplus_{i \in I_0} \Z c_i \varpi_i.
\]
Let $T(M)$ and $T(\wti{M})$ be the subgroups of $\mathrm{GL}(\fh^*)$ defined by 
\[ T(M) = \{ t_{\gl} \mid \gl \in M\}, \ \ \ T(\wti{M}) = \{ t_{\gl} \mid \gl \in \wti{M} \} . 
\]
It is known that \cite[Proposition 6.5]{MR1104219}
\[ W \cong W_0 \ltimes T(M).
\]
Define the subgroup $\widetilde{W}$ of $\mathrm{GL}(\fh^*)$ by
\[ \widetilde{W} = W_0 \ltimes T(\widetilde{M}),
\]
which is called the \textit{extended affine Weyl group}.
The action of $\wti{W}$ preserves $\gD$, and elements $w \in W_0$ and $\gl \in \widetilde{M}$ satisfy
\[ wt_{\gl} w^{-1} = t_{w(\gl)}.
\]
In the sequel, we assume that the positive integer $N$ in (\ref{eq:def_of_P}) satisfies
\begin{equation*}\label{eq:preserve}
  2^{-1}\big(\aff (\gl) , \aff (\gl)\big) \in N^{-1} \Z \ \ \ \text{for all} \ \gl \in \wti{M},
\end{equation*}
which ensures that $\wti{W}$ preserves $P$.

Let $\Aut (\gG)$ be the group of automorphisms of the Dynkin diagram $\gG$, that is, the group of permutations $\tau$ of $I$ 
satisfying $a_{ij} = a_{\tau(i) \tau(j)}$ for all $i,j \in I$.
Note that $\tau \in \Aut(\gG)$ satisfies
\begin{equation*} \label{eq:property_of_a}
  a_{\tau(i)} = a_i \ \text{and} \ a_{\tau(i)}^{\vee} = a_i^{\vee} \ \ \  \text{for all} \ i \in I.
\end{equation*}
Let  $C \subseteq \fh^*_{\R} = \R \otimes_\Z P$ be the fundamental chamber
(i.e.\ $C = \{ \gl \in \fh^*_\R \mid ( \gl, \ga_i )  \ge 0 \ \text{for all} \ i \in I \}$),
and $\gS \subseteq \wti{W}$ the subgroup consisting of elements preserving $C$.
Then we have
\[ \wti{W} \cong W \rtimes \gS.
\]
The length function $\ell$ is exnteded on $\wti{W}$ by setting $\ell(w\tau) = \ell(w)$ for $w \in W$ and $\tau \in \gS$.
Note that an element $w \in \wti{W}$ belongs to $\gS$ if and only if $w$ preserves the set of simple roots 
$\{\ga_0,\dots,\ga_n\}$.
Hence $\tau \in \gS$ induces a permutation of $I$ (also denoted by $\tau$) by $\tau(\ga_i) = \ga_{\tau(i)}$,
which belongs to $\Aut(\gG)$ since $( \ , \ )$ is $\tau$-invariant.
By abuse of notation, we denote by $\gS$ both the subgroup of $\wti{W}$ and the subgroup of $\Aut(\gG)$.

We shall describe the subgroup $\gS \subseteq \Aut (\gG)$ explicitly.
A node $i \in I$ is called a \textit{special} node if $i \in \Aut(\gG)\cdot 0$.
Let $I^s \subseteq I$ be the set of special nodes.
$I^s$ for nonexceptional $\fg$ are as follows:
\[ I^s = \begin{cases} \{0,1,\dots,n\} & \text{for} \ A_n^{(1)}, \\
                 \{0,1\}       & \text{for} \ B_n^{(1)}, \ A_{2n-1}^{(2)}, \\
                 \{0,n\}       & \text{for} \ C_n^{(1)},  \ D_{n+1}^{(2)}, \\
                 \{0,1,n-1,n\} & \text{for} \ D_n^{(1)}, \\
                 \{0\}          & \text{for} \ A_{2n}^{(2)}.
   \end{cases}
\]
Assume that $i \in I^s\setminus \{ 0\}$ (in particular $\fg \neq A_{2n}^{(2)}$), and define $\tau^i \in \wti{W}$ by 
\[ \tau^i = t_{\varpi_i}w_i,
\]
where $w_i$ denotes the unique element of $W_0$ which maps the simple system $\{\ga_1, \ga_2, \dots, \ga_n\}$ of $\fg_0$ to
$\{-\gt, \ga_1, \dots, \Hat{\ga_i}, \dots, \ga_n\}$ with $\gt = \gd - \ga_0 \in \gD$.
We put $\tau^0 =\id$.
The following proposition is well-known, but  we give the proof for completeness:

\begin{Prop}\label{Prop:special_auto}
  {\normalfont(i)} For all $i \in I^s$, $\tau^i$ belongs to $\gS$. \\
  {\normalfont(ii)} The map $I^s \to \gS$ defined by $i \mapsto \tau^i$ is bijective. \\
  {\normalfont(iii)} If $\tau\in\gS$ satisfies $\tau(i) = 0$, then we have $\tau= (\tau^i)^{-1}$.
\end{Prop}

\begin{proof}
  If $\fg$ is of type $A_{2n}^{(2)}$, then $I^s = \{0\}$ and $\wti{M} = M$, which obviously imply the assertions.
  So we may assume that $\fg$ is not of this type.
  (i) 
  Let $i \in I^s \setminus \{ 0\}$,
  and recall that $w_i$ maps $\{\ga_1, \ga_2, \dots, \ga_n\}$ to $\{-\gt, \ga_1, \dots, \Hat{\ga_i}, \dots, \ga_n\}$.
  Then $w_i(-\gt) = \ga_i$ also holds,
  and hence it is easily checked from the equation (\ref{eq:translation}) that $\tau^i = t_{\varpi_i}w_i$ preserves the set 
  $\{\ga_0,\dots, \ga_n\}$, which implies $\tau^i$ belongs to $\gS$.
  (ii) The injectivity is obvious.
  Let $\tau \in \gS \setminus \{ \id \}$ be an arbitrary element, 
  and decompose it as $\tau = t_{\gl_{\tau}}w_{\tau}$ where $\gl_{\tau} \in \wti{M}$
  and $w_{\tau} \in W_0$.
  Since $t_{\gl_{\tau}}$ acts trivially on $P_{\cl}^0$, we have $w_{\tau}\big(\cl(\ga_{j})\big) 
  = \cl\big(\ga_{\tau(j)}\big)$ for $j \in I$,
  which implies $w_\tau=w_{\tau(0)}$.
  Then since 
  \[ t_{\gl_{\tau}}(\ga_{j}) =  \tau w_{\tau(0)}^{-1}(\ga_j) 
     = \tau(\ga_{\tau^{-1}(j)} - \gd_{j,\tau(0)}\gd) = \ga_j - \gd_{j,\tau(0)}\gd \ \ \ \text{for} \ j \in I_0,
  \]
  (\ref{eq:translation}) forces $\gl_{\tau} = \varpi_{\tau(0)}$, and the surjectivity follows.
  From the proof of (ii), we see that $\tau(0) = i$ implies $\tau = \tau^i$.
  Hence the assertion (iii) follows.
\end{proof}

For nonexceptional $\fg$, $\tau^i$ for $i \in I^s \setminus \{0\}$  are as follows:
\begin{enumerate}
  \setlength{\parskip}{0pt} 
	\setlength{\itemsep}{1pt} 
  \item[$A_n^{(1)}$:] $\tau^i(j) = j +i$ mod $n+1$ for all $j \in I$.
  \item[$B_n^{(1)}$,] $D_{n+1}^{(2)}$: $\tau^1= (0,1)$.
  \item[$C_n^{(1)}$,] $A_{2n-1}^{(2)}$: $\tau^n(j) = n-j$ for all $j \in I$.
  \item[$D_n^{(1)}$,]  $n$ odd: $\tau^1 = (0,1)(n-1,n)$, \\
  $\tau^{n-1}(0,1,n-1,n) = (n-1,n,1,0)$, $\tau^{n-1}(j) = n-j$ for $j \in I \setminus I^s$, \\
  $\tau^n(0,1,n-1,n) = (n,n-1,0,1)$, $\tau^n(j) = n-j$ for $j \in I \setminus I^s$.
  \item[$D_n^{(1)}$,] $n$ even: $\tau^1 = (0,1)(n-1,n)$, \\
  $\tau^{n-1}(0,1,n-1,n) = (n-1,n,0,1)$, $\tau^{n-1}(j) = n-j$ for $j \in I \setminus I^s$, \\
  $\tau^n(0,1,n-1,n) = (n,n-1,1,0)$, $\tau^{n}(j) = n-j$ for $j \in I \setminus I^s$.
\end{enumerate}

In the sequel, we assume that the fundamental weights $\gL_0, \dots, \gL_n$ are chosen to satisfy $\tau(\gL_j) = \gL_{\tau(j)}$
for all $\tau \in \gS$ and $j \in I$.
This is always possible by choosing $\gL_j$ arbitrarily for representatives of $I/ \gS$ 
and setting $\gL_{\tau(j)} = \tau(\gL_j)$ for $\tau \in \gS$.
Then each element $\tau' \in \Aut(\gG)$ acts on $P$ by $\tau'(\gL_i) = \gL_{\tau'(i)}$ and $\tau'(\gd) = \gd$.

\section{Definition of crystals}\label{Section:crystals}

Let $U_q(\fg)$ be the quantum affine algebra associated with $\fg$ and $U_q'(\fg)$ the one without the degree operator.
The weight lattices of $U_q(\fg)$ and $U_q'(\fg)$ are $P$ and $P_{\cl}$ respectively.

A $U_q(\fg)$-crystal (resp.\ $U_q'(\fg)$-crystal) is by definition a set $B$ 
equipped with weight function $\wt\colon B \to P$ (resp.\ $\wt\colon B \to P_{\cl}$) and
Kashiwara operators $e_i,f_i\colon B \to B \sqcup \{0\}$ for $i \in I$  satisfying
\begin{align*}
  \wt(e_i b) &= \wt(b) + \ga_i \ \text{and} \ f_i(e_i b) = b \ \ \ \text{for all} \ i \in I, b \in B \ \text{such that}
  \ e_i b \neq 0, \\
  \wt(f_i b) &= \wt(b) - \ga_i \ \text{and} \ e_i(f_i b) = b \ \ \ \text{for all} \ i \in I, b \in B \ \text{such that} 
  \ f_i b \neq 0,
\end{align*}
and $\langle \wt(b), \ga_i^{\vee} \rangle = \gph_i(b)- \gee_i(b)$ where 
\[ \gee_i(b) = \max\{ k \ge 0 \mid e_i^k b \neq 0 \}, \ \ \ \gph_i(b) = \max\{ k \ge 0 \mid f_i^k b \neq 0 \}.
\]
In this article, we always assume that $\gee_i(b) < \infty$ and $\gph_i(b) < \infty$. 
We call $B$ a crystal if $B$ is either a $U_q(\fg)$-crystal or a $U_q'(\fg)$-crystal.

\begin{Rem} \label{Rem:restriction} \normalfont
  A $U_q(\fg)$-crystal $B$ can be regarded naturally as a $U_q'(\fg)$-crystal by replacing the weight function 
  $\wt\colon B \to P$ by $\cl \circ \wt\colon B \to P_{\cl}$.
\end{Rem}

For two crystals $B_1$ and $B_2$, their tensor product $B_1 \otimes B_2= \{ b_1 \otimes b_2 \mid b_1 \in B_1, b_2 \in B_2\}$
is defined with weight function $\wt(b_1 \otimes b_2) = \wt(b_1) + \wt(b_2)$ and Kashiwara operators
\begin{align*}\label{eq:tensor}
   e_i(b_1 \otimes b_2) &= \begin{cases} e_i b_1 \otimes b_2 & \text{if} \ \varphi_i(b_1) \ge \gee_i(b_2), \\
                                        b_1 \otimes e_i b_2 & \text{if} \ \varphi_i(b_1) < \gee_i(b_2),
                          \end{cases} \\
   f_i(b_1 \otimes b_2) &= \begin{cases} f_i b_1 \otimes b_2 & \text{if} \ \varphi_i(b_1) > \gee_i(b_2), \\
                                        b_1 \otimes f_i b_2 & \text{if} \ \varphi_i(b_1) \le \gee_i(b_2).
                          \end{cases}
\end{align*}
Note that it follows for $i \in I$ that
\begin{align}
   \gee_i(b_1 \otimes b_2) &= \gee_i(b_1) + \max\big\{0, \gee_i(b_2) - \gph_i(b_1)\big\}, \\
   \gph_i(b_1 \otimes b_2) &= \gph_i(b_2) + \max\big\{0, \gph_i(b_1) - \gee_i(b_2) \big\}. \nonumber
\end{align}
The following lemma is obvious:

\begin{Lem}\label{Lem:add}
  Let $B_1, B_2$ be crystals, $b_j \in B_j$ {\normalfont(}$j =1,2${\normalfont)} and $i \in I$.
  Put $m = \max\{0, \gee_i(b_2) - \varphi_i(b_1)\}$. Then we have
  \[ e_i^{m+1}(b_1 \otimes b_2) = e_ib_1 \otimes e_i^m b_2.
  \]
\end{Lem}

For crystals $B_1, B_2$ and their subsets $S_j \subseteq B_j$,
a bijection $\Psi\colon S_1 \to S_2$ is said to be an \textit{isomorphism of full subgraphs}
if $\wt \Psi(b) = \wt(b)$ for $b\in S_1$, $\Psi(e_ib) = e_i \Psi(b)$ for $b \in S_1$ such that 
$e_ib \in S_1 \sqcup \{ 0 \}$, and $\Psi(f_ib) = f_i \Psi(b)$ for $b \in S_1$ such that $f_ib \in S_1 \sqcup \{ 0 \}$.
Here $\Psi(0)$ is understood as $0$. If there exists an isomorphism $S_1 \to S_2$ of full subgraphs, 
we say $S_1$ and $S_2$ are \textit{isomorphic as full subgraphs} and write $S_1 \cong S_2$.

For a crystal $B$ and $\tau \in \Aut (\gG)$, we define a crystal $\widetilde{\tau}(B)$ as follows: 
as a set $\widetilde{\tau}(B) = \{\widetilde{\tau}(b) \mid b \in B \}$, where $\wti{\tau}(b)$ is just a symbol.
Its weight function and Kashiwara operators are defined by 
\begin{equation}\label{eq:weight_change}
  \wt\big(\widetilde{\tau}(b)\big) = \tau \big(\wt(b)\big) \ \text{and}
\end{equation}
\[ e_i \widetilde{\tau}(b) = \widetilde{\tau}(e_{\tau^{-1}(i)}b), \ \ \ f_i \widetilde{\tau}(b) 
   = \widetilde{\tau}(f_{\tau^{-1}(i)}b),
\]
where $\widetilde{\tau}(0)$ is understood as $0$.
Obviously we have
\[ \widetilde{\tau}(B_1 \otimes B_2) \cong \widetilde{\tau}(B_1) \otimes \widetilde{\tau}(B_2)
\]
for two crystals $B_1$ and $B_2$.
For a subset $S \subseteq B$, a subset $\widetilde{\tau}(S) \subseteq \widetilde{\tau}(B)$
is defined in the obvious way.

For $J \subseteq I$, we denote by $U_q(\fg_J)$ the subalgebra of $U_q(\fg)$ whose simple roots are $J$.
If $J = I_0$, we denote $U_q(\fg_J)$ by $U_q(\fg_0)$.
$U_q(\fg_J)$-crystals are defined in a similar way.
For a crystal $B$ and a proper subset $J$ of $I$, a connected component of $B$ regarded as a $U_q(\fg_J)$-crystal
is called a $U_q(\fg_J)$-\textit{component} of $B$.

\begin{Def}[\cite{MR1607008}] \normalfont 
  We say a crystal $B$ is \textit{regular} 
  if for every proper subset $J$ of $I$, $B$ regarded as a $U_q(\fg_J)$-crystal is 
  isomorphic to a direct sum of the crystal bases of integrable highest weight $U_q(\fg_J)$-modules.
\end{Def}

Let $J \subseteq I$. 
For a crystal $B$, we say that $b \in B$ is \textit{$U_q(\fg_J)$-highest weight} if $e_j b = 0$ for all $j \in J$.
For a proper subset $J$ of $I$ and a regular crystal $B$, 
every $U_q(\fg_J)$-component of $B$ contains a unique $U_q(\fg_J)$-highest weight element.

By \cite{MR1262212}, the actions of simple reflections on a regular crystal $B$ defined by
\[ S_{s_i}(b) = \begin{cases} f_i^{\langle \wt(b), \ga_i^{\vee} \rangle} b & \text{if} \ \langle \wt(b), \ga_i^{\vee} \rangle 
                              \ge 0, \\
                              e_i^{-\langle \wt(b), \ga_i^{\vee} \rangle} b & \text{if} \ \langle \wt(b), \ga_i^{\vee} \rangle
                              < 0
                \end{cases}
\]
are extended to the action of $W$ denoted by $w \mapsto S_w$.
For every $w \in W$ and $b\in B$, we have $\wt\big(S_w(b)\big) = w\big(\wt(b)\big)$.
An element $b \in B$ is called \textit{extremal} if 
for every $w \in W$ and $i \in I$, 
\[ e_iS_w(b) = 0 \ \text{if} \ \left\langle \wt\big(S_w(b)\big), \ga_i^{\vee} \right\rangle \ge 0 \ \text{and} \
   f_iS_w(b) = 0 \ \text{if} \ \left\langle \wt\big(S_w(b)\big), \ga_i^{\vee} \right\rangle \le 0.
\]

\section{Demazure crystals}\label{section:Demazure_crystals}

For a subset $S$ of a crystal $B$ and $i \in I$, we denote 
$\mathcal{F}_iS = \{ f_i^k b \mid b \in S, k \ge 0 \} \setminus \{ 0 \} \subseteq B$.

For $\gL \in P^+$, let $V(\gL)$ denote the integrable highest weight $U_q(\fg)$-module with highest weight $\gL$,
and $B(\gL)$ its crystal basis with highest weight element $u_{\gL}$.
Let $w$ be an element of $W$ and $w = s_{i_k} s_{i_{k-1}}\cdots s_{i_1}$ its reduced expression.
Then it is known that the subset
\[ B_w(\gL) = \cF_{i_k} \cF_{i_{k-1}} \cdots \cF_{i_1} \{u_{\gL}\}  \subseteq B(\gL)
\]
is independent of the choice of the reduced expression of $w$ \cite{MR1240605}.

\begin{Def} \normalfont
  The subset $B_w(\gL)$ of $B(\gL)$ is called the \textit{Demazure crystal} associated with $\gL$ and $w$.
\end{Def}

\begin{Rem} \normalfont \label{Rem!}
  Let $\fb$ be the standard Borel subalgebra of $\fg$ and $U_q(\fb) \subseteq U_q(\fg)$ 
  the corresponding quantized enveloping algebra.
  The \textit{Demazure module} $V_w(\gL)$ is defined by the $U_q(\fb)$-submodule of $V(\gL)$ generated by 
  the weight space $V(\gL)_{w(\gL)}$.
  The Demazure crystal $B_w(\gL)$ is known to be the crystal basis of $V_w(\gL)$
  in a suitable sense \cite{MR1240605}, which is why it is so named.
\end{Rem}

For a subset $S$ of a crystal and $w \in W$ with a reduced expression $w = s_{i_k}\cdots s_{i_1}$,
we write $\cF_w S = \cF_{i_k} \cdots \cF_{i_1} S$ if it is well-defined.
For example, $\cF_w \{ u_{\gL}\} = B_w(\gL)$.

\begin{Lem}\label{Lem:Demazure}
  Let $\gL \in P^+$ and $w \in W$. \\
  {\normalfont(i)} We have $\widetilde{\tau} B_w(\gL) \cong B_{\tau w \tau^{-1}} \big(\tau (\gL)\big)$ 
                   for $\tau \in \Aut(\gG)$. \\
  {\normalfont(ii)} For $i \in I$, we have
  \begin{equation}\label{eq;Demazure}
    \cF_i B_w(\gL) = \begin{cases} B_w(\gL) & \text{if} \ \ell(s_iw) = \ell(w) - 1, \\
                                    B_{s_iw}(\gL) & \text{if}\ \ell(s_i w) = \ell(w) +1.
                     \end{cases}
  \end{equation}
  {\normalfont(iii)} For every $w' \in W$, $\cF_{w'} B_w(\gL)$ is well-defined, and
                     $\cF_{w'} B_w(\gL) = B_{w''}(\gL)$ for some $w'' \in W$.
                     Moreover, if $\ell(w'w) = \ell(w') + \ell(w)$, then $w'' = w'w$.
\end{Lem}

\begin{proof}
  Since $\widetilde{\tau}B(\gL) \cong B(\tau(\gL))$ and $\wti{\tau}( \cF_i S) = \cF_{\tau(i)} \tilde{\tau}(S)$
  for every $S \subseteq B(\gL)$, (i) follows.
  When $\ell(s_iw) = \ell(w) + 1$, (\ref{eq;Demazure}) follows by definition,
  and when $\ell(s_iw) = \ell(w)-1$, this follows since
  \[ \cF_iB_w(\gL) = \cF_i\big(\cF_iB_{s_iw}(\gL)\big) = \cF_iB_{s_iw}(\gL) = B_w(\gL).
  \]
  The assertion (ii) is proved.
  To see that $\cF_{w'} B_w(\gL)$ is well-defined, it suffices to show 
  the operators $\cF_i$ on Demazure crystals satisfy braid relations: if the order of $s_i s_j$ for $i,j \in I$ $(i\neq j)$
  is $m < \infty$, then we have $\underbrace{\cF_i\cF_j\cF_i\cdots}_m B_w(\gL) = 
  \underbrace{\cF_j\cF_i\cF_j\cdots}_m B_w(\gL)$.
  Since the element $\underbrace{s_is_js_i\cdots}_m = \underbrace{s_js_is_j\cdots}_m$ is the longest element of
  the subgroup $W_{i,j} = \langle s_i, s_j \rangle \subseteq W$, (ii)  
  implies 
  \[ \underbrace{\cF_i\cF_j\cF_i\cdots}_m B_w(\gL) = B_{w''}(\gL)= \underbrace{\cF_j\cF_i\cF_j\cdots}_m B_w(\gL),
  \] 
  where $w''$ is the unique element of the set $\{ \gs w \mid \gs \in W_{i,j} \}$ whose length is maximal. 
  Hence our assertion is proved. 
  The remaining statements of (iii) are obvious from (ii). 
\end{proof}

For $w \in W$ and $\tau \in \Aut(\gG)$, we write $\cF_{w\tau} = \cF_{w} \widetilde{\tau}$
and $B_{w\tau}(\gL) = B_w\big(\tau(\gL)\big)$ for the notational convenience.
The following proposition is immediate from Lemma \ref{Lem:Demazure}.

\begin{Prop}\label{Prop:twisted_Demazure}
  For every $\gL \in P^+$ and $w,w' \in \wti{W}$, there exists $w'' \in \wti{W}$ such that
  \[ \cF_{w'}B_w(\gL) \cong B_{w''}(\gL).
  \]
  Moreover, if $\ell(w'w) = \ell(w')+ \ell(w)$, then $w'' = w'w$.
\end{Prop}

Let $\C[P]$ denote the group algebra of $P$ with basis $e^\gl$ ($\gl \in P$),
and define for $i \in I$ a linear operator $D_i$ on $\C[P]$ by
\begin{equation*}
   D_i(f) = \frac{f - e^{-\ga_i}\cdot s_i(f)}{1-e^{-\ga_i}},
\end{equation*}
where $s_i$ acts on $\C[P]$ by $s_i(e^{\gl}) = e^{s_i(\gl)}$.
The operator $D_i$ is called the \textit{Demazure operator} associated with $i$.
Note that $D_i(f) = f$ holds if $f$ is $s_i$-invariant.
From this, it is easily checked that $D_i^2 = D_i$.

For every reduced expression $w = s_{i_k} \cdots s_{i_1}$ of $w \in W$, 
the operator $D_w = D_{i_k} \cdots D_{i_1}$ on $\C[P]$ is independent
of the choice of the expression \cite{MR1923198}.
The weight sum of a Demazure crystal is known to be expressed using Demazure operators:

\begin{Thm}[\cite{MR1240605}] \label{Thm:Kashiwara}
  For $\gL \in P^+$ and $w \in W$, we have
  \[ \sum_{b \in B_w(\gL)} e^{\wt(b)}= D_w (e^\gL).
  \]
\end{Thm}

For $w \in W$ and $\tau \in \Aut(\gG)$, we define an operator $D_{w\tau}$ on $\C[P]$ by $D_{w\tau} = D_w\circ \tau$,
where $\tau$ acts on $\C[P]$ by $\tau(e^{\gl}) = e^{\tau(\gl)}$.

\begin{Cor} \label{Cor:Demazure_character}
  Let $S$ be a disjoint union of Demazure crystals and $i \in I$.
  For every $w \in \wti{W}$ we have
  \begin{equation} \label{eq:Demazure_character} 
    \sum_{b \in \cF_w(S)} e^{\wt(b)} = D_w \left( \sum_{b \in S} e^{\wt(b)} \right).
  \end{equation}
\end{Cor}

\begin{proof}
  We may assume that $S$ is a single Demazure crystal, say $S = B_{w'}(\gL)$.
  By Proposition \ref{Prop:twisted_Demazure}, 
  it suffices to show the assertion for $w = \tau \in \gS$ and $w = s_i$ for $i \in I$.
  When $w = \tau$, the assertion is obvious from (\ref{eq:weight_change}).
  Assume that $w = s_i$.
  If $\ell (s_i w') = \ell(w') + 1$, then we have $\cF_iB_{w'}(\gL) = B_{s_iw'}(\gL)$, and 
  the assertion follows from Theorem \ref{Thm:Kashiwara}.
  If $\ell(s_i w') = \ell(w') -1$, then we have $\cF_iB_{w'}(\gL) = B_{w'}(\gL)$.
  On the other hand, it follows that
  \[ D_i \left( \sum_{b \in B_{w'}(\gL)} e^{\wt(b)} \right) = \sum_{b \in B_{w'}(\gL)} e^{\wt(b)}
  \]
  since the weight sum
  \[ \sum_{b \in B_{w'}(\gL)} e^{\wt(b)} = D_i\left(\sum_{b \in B_{s_iw'}(\gL)} e^{\wt(b)}\right)
  \]
  is $s_i$-invariant.
  Hence the assertion follows.
\end{proof}

It is known that $B(\gL)\otimes B(\gL')$ for $\gL, \gL' \in P^+$ is isomorphic to 
a direct sum of the crystal bases of integrable highest weight 
modules, that is,
\begin{equation}\label{eq:isom}
  B(\gL) \otimes B(\gL') \cong \bigoplus_{\gl \in T} B(\gl),
\end{equation}
where $T$ is a possibly infinite multiset of elements of $P^+$. 
The following theorem, which was proved in \cite[Proposition 12]{MR1887117}
and \cite[Theorem 2.11]{MR1987017}, is known as the combinatorial excellent filtration theorem:

\begin{Thm} \label{Thm:excellent_filtration}
  The image of the subset $u_{\gL} \otimes B_w(\gL')$ of $B(\gL) \otimes B(\gL')$ under the isomorphism 
  {\normalfont(\ref{eq:isom})}
  is a disjoint union of Demazure crystals.
\end{Thm}

For $\gL \in P^+$, $B(\gL)$ is regular and $u_\gL$ is extremal.
For $w \in W$ and $\tau \in \Aut(\gG)$, 
set $u_{w\tau(\gL)} = S_w\big(u_{\tau(\gL)}\big) \in B\big(\tau(\gL)\big)$.
By definition we have $u_{w\tau(\gL)} \in B_{w\tau}(\gL)$.
Later we need the following lemma:

\begin{Lem}\label{Lem}
  Let $\gL \in P^+$ and $w \in \wti{W}$, and assume that $\langle w(\gL), \ga_i^\vee\rangle \le 0$ for all $i \in I_0$.
  Then for any $b \in B_w(\gL)$, we have $b = u_{w(\gL)}$ or 
  \[ \cl\big(\wt(b)\big) \in \cl\big(w(\gL)\big) + \Big(Q_0^+ \setminus \{0\}\Big).
  \]
\end{Lem}

\begin{proof}
  Assume that $w = w'\tau$ with $w' \in W$ and $\tau \in \gS$.
  In view of Remark \ref{Rem!}, it follows that $\wt(b) \in w(\gL) + Q^+$.
  Then since $\wt(b)$ is a weight of $V\big(\tau(\gL)\big)$, it is proved from \cite[Proposition 11.3 (a)]{MR1104219} that 
  $\wt(b) = w(\gL)$, or $\wt(b) \in w(\gL) + \Big(Q_0^+ \setminus \{0\}\Big) + \Z_{\ge 0}\gd$.
  Hence the assertion follows.
\end{proof}

\section{Perfect Kirillov-Reshetikhin crystals}

\textit{From this section to the end of the article, we assume that the type of $\fg$ is nonexceptional}
(i.e.\ one of the types
$A_n^{(1)}, B_n^{(1)}, C_n^{(1)}, D_n^{(1)}, A_{2n-1}^{(2)}, A_{2n}^{(2)}, D_{n+1}^{(2)}$).
Note that some of the statements below on Kirillov-Reshetikhin crystals may
have not been proved or not be true for exceptional $\fg$.

\subsection{Kirillov-Reshetikhin crystals}

For a $U_q'(\fg)$-crystal $B$, define two maps $\gee, \gph\colon B \to P_{\cl}^+$ by
\[ \gee(b) = \sum_{i \in I} \gee_i(b) \gL_i, \ \ \ \gph(b) = \sum_{i \in I} \gph_i(b) \gL_i \ \ \ \text{for} \ b \in B.
\]
Note that $\wt(b) = \gph(b) - \gee(b)$.

Kirillov-Reshetikhin modules $W^{r,\ell}$ (KR modules for short) are irreducible finite-dimensional $U_q'(\fg)$-modules 
parametrized by $r \in I_0$ and $\ell \in \Z_{\ge 1}$ (see \cite{MR1903978} for the precise definition).
For nonexceptional $\fg$, the following theorem is known:

\begin{Thm}[\cite{MR1194953,MR2389790,MR2403558}] \label{Thm:KR}
  For each $r \in I_0$ and $\ell \in \Z_{\ge 1}$, the KR module $W^{r,\ell}$ has a crystal basis $B^{r,\ell}$.
\end{Thm}

The crystals $B^{r,\ell}$ are called the \textit{Kirillov-Reshetikhin crystals} (KR crystals for short).
In this article we denote by $\mathcal{C}$ the set consisting of tensor products of KR crystals.

\begin{Def}[\cite{MR1607008}]  \normalfont
  A finite regular $U_q'(\fg)$-crystal $B$ is called \textit{simple} if there exists $\gl \in P_{\cl}^0$ such that 
  $B$ has a unique element whose weight is $\gl$, the weights of
  $B$ are contained in the convex hull of $W \gl$, and the weight of each extremal element is in $W\gl$.
\end{Def}

\begin{Prop}[{\cite[Proposition.3.8 (1)]{MR2945586}}] \label{Prop:simplicity}
  Every $B \in \cC$ is simple. 
\end{Prop}

Since $B \in \cC$ is simple, $B$ has a unique extremal element $u(B)$ such that $\wt \big(u(B)\big) \in - P_0^+$.
It is known that $u(B^{r,\ell})$ is the unique element with weight $\ell w_0(\varpi_r)$, 
and we have $u(B_1 \otimes B_2) = u(B_1) \otimes u(B_2)$ for $B_1, B_2 \in \cC$.
Every $B \in \cC$ is connected by \cite[Lemma 1.9 and 1.10]{MR1607008}.
Then by \cite{MR1641035}, we have the following:

\begin{Lem}[{\cite[Lemma 3.3 (b)]{MR1641035}}] \label{Lem:strongly_connected} \normalfont
  For $B \in \cC$ and every $b \in B$, we have
  \[ B = \{ e_{i_k} \cdots e_{i_1} (b) \mid k \ge 0, i_j \in I \} \setminus \{ 0 \}.
  \]
\end{Lem}

The following proposition is important:

\begin{Prop} \label{Prop:ST} 
  Let $B \in \cC$. For every $\tau \in \gS$, there exists a unique isomorphism
  $\rho_{\tau}\colon \widetilde{\tau}(B) \stackrel{\sim}{\to} B$
  of $U'_q(\fg)$-crystals.
\end{Prop}

\begin{proof}
  For a single KR crystal $B = B^{r,\ell}$, $\wti{\tau}(B^{r,\ell}) \cong B^{r,\ell}$ was proved in \cite[Lemma 6.5]{ST}.
  This implies $\wti{\tau}(B) \cong B$ for general $B \in \cC$ 
  since $\wti{\tau}(B_1\otimes B_2) \cong \wti{\tau}(B_1) \otimes \wti{\tau}(B_2)$.
  Since $B$ is connected and an element of $B$ with weight $\wt\big(u(B)\big)$ is unique,
  the uniqueness of the isomorphism holds.
\end{proof}

Using the isomorphism $\rho_\tau$ in the above proposition, 
we define an action of $\gS$ on $B \in \cC$ by $\tau(b) = \rho_\tau \big( \wti{\tau}(b) \big)$ 
for $\tau \in \gS$.
This action satisfies 
\begin{equation}\label{eq:tau}
   \tau \circ e_i = e_{\tau(i)} \circ \tau \ \text{and} 
   \ \tau \circ f_i = f_{\tau(i)} \circ \tau \ \ \ \text{for all} \ i \in I.
\end{equation}

\begin{Lem}\label{Lem:highest}
  For every $\tau \in \gS$, there exists some $w \in W_0$ such that 
  \[ \tau\big(u(B)\big) = S_w \big(u(B)\big) \ \ \ \text{for all} \ B \in \cC.
  \]
\end{Lem}

\begin{proof}
  Since $\tau \in \wti{W} = W_0 \ltimes T(\wti{M})$ and $T(\wti{M})$ acts on $P_{\cl}^0$ trivially,
  there exists $w \in W_0$ such that $\tau|_{P_{\cl}^0} = w|_{P_{\cl}^0}$.
  Then since
  \[ \wt\Big(\tau\big(u(B)\big)\Big) = \tau\Big( \wt\big( u(B)\big)\Big) = w\Big(\wt\big(u(B)\big)\Big) = \wt\Big(S_w\big(u(B)\big)\Big),
  \]
  $\tau\big(u(B)\big) = S_w \big(u(B)\big)$ follows by Proposition \ref{Prop:simplicity}.
\end{proof}

The $U_q(\fg_0)$-crystal structure of a KR crystal is known by \cite{MR1836791,MR2576287}.
In particular, we have the following proposition (for nonexceptional $\fg$):

\begin{Prop}
  A KR crystal $B^{r,\ell}$ is multiplicity free as a $U_q(\fg_0)$-crystal.
  In other words, any two distinct $U_q(\fg_0)$-components of $B^{r,\ell}$ are not isomorphic as $U_q(\fg_0)$-crystals.
\end{Prop}

\begin{Cor}\label{Cor:multiplicity_free}
  Let $b_1, b_2 \in B^{r,\ell}$ be two distinct $U_q(\fg_0)$-highest weight elements.
  Then  we have 
  \[ \gph(b_1) - \gph(b_2) \notin \Z \gL_0.
  \]
\end{Cor}

\begin{proof}
  For $j = 1,2$, let $B_j \subseteq B^{r,\ell}$ be the $U_q(\fg_0)$-component containing $b_j$.
  Then as a $U_q(\fg_0)$-crystal, $B_j$ is isomorphic to the crystal basis of the integrable highest weight 
  $U_q(\fg_0)$-module with highest weight $\sum_{i \in I_0} \gph_i(b_j) \varpi_i$.
  Now the assertion is obvious from the above proposition.
\end{proof}

\subsection{Perfect KR crystals}\label{subsection:Perfect_KR_crystals}

For a $U_q'(\fg)$-crystal $B$ such that $\wt(B) \subseteq P_{\cl}^0$, we define the \textit{level} of $B$ by
\[ \lev (B) = \min_{b \in B} \langle \gph(b), K \rangle =\min_{b \in B} \langle \gee(b), K \rangle,
\]
and the subset $B_{\min}$ by
\begin{align*} B_{\min} &= \{ b \in B \mid \langle \gph(b), K \rangle = \lev (B)\} \\
                       &= \{ b \in B \mid \langle \gee(b), K \rangle = \lev (B)\}.
\end{align*}

\begin{Def}[\cite{MR1187560}] \normalfont
  For a positive integer $\ell$, a $U_q'(\fg)$-crystal $B$ is called a \textit{perfect crystal of level} $\ell$
  if $B$ satisfies the following conditions: \\[-0.6cm]
  \begin{enumerate}
    \item[(i)] $B$ is isomorphic to the crystal basis of a finite-dimensional $U_q'(\fg)$-module. \\[-0.6cm]
    \item[(ii)] $B \otimes B$ is connected.\\[-0.6cm]
    \item[(iii)] There exists $\gl \in P_{\cl}^0$ such that $\wt(B) \subseteq \gl - \sum_{i \in I_0} \Z_{\ge 0} \ga_i$
                 and there exists a unique element in $B$ with weight $\gl$. \\[-0.6cm]
    \item[(iv)] The level of $B$ is $\ell$. \\[-0.6cm]
    \item[(v)]  Both the maps $\gee$ and $\gph$ induce bijections between the set $B_{\min}$ and $(P_\cl^+)^\ell$.
  \end{enumerate}     
\end{Def}

The following lemma is immediate:

\begin{Lem} \label{Lem:level}
  Let $B_1, B_2$ be perfect crystals.\\[0.2cm]
  {\normalfont(i)} $\lev(B_1 \otimes B_2) = \max\{\lev(B_1), \lev(B_2)\}$.  \\[0.1cm]
  {\normalfont(ii)} If $\lev(B_1) \ge \lev(B_2)$, then $b_1 \otimes b_2 \in B_1 \otimes B_2$ belongs 
                    to $(B_1\otimes B_2)_{\min}$ if and only if $b_1 \in (B_1)_{\min}$ and $\gph(b_1) - \gee(b_2) \in P_\cl^+$.
                    Moreover if $b_1 \otimes b_2 \in (B_1\otimes B_2)_{\min}$, then 
                    \[ \gee(b_1 \otimes b_2) = \gee(b_1), \ \ \ \gph(b_1 \otimes b_2) = \gph(b_1) + \wt(b_2).
                    \]
  {\normalfont(iii)} If $\lev(B_1) \le \lev(B_2)$, then $b_1 \otimes b_2 \in B_1 \otimes B_2$ belongs
                    to $(B_1\otimes B_2)_{\min}$ if and only if $b_2 \in (B_2)_{\min}$ and $\gee(b_2) - \gph(b_1) \in P_\cl^+$.
                    Moreover if $b_1 \otimes b_2 \in (B_1\otimes B_2)_{\min}$, then 
                    \[ \gee(b_1 \otimes b_2) = \gee(b_2) - \wt(b_1), \ \ \ \gph(b_1 \otimes b_2) = \gph(b_2).
                    \]
\end{Lem}

The significance of the perfectness is due to the following theorem:

\begin{Thm}[\cite{MR1187560}] \label{Thm:Path}
  Let $B$ be a perfect crystal of level $\ell$, $\gL\in P^+$ a dominant integral weight of level $\ell$,
  and $b$ the unique element of $B$ satisfying $\gee(b) = \cl(\gL)$.
  Then for all $\gL' \in P^+$ such that $\gph(b) = \cl(\gL')$, we have
  \[ B(\gL) \otimes B \stackrel{\sim}{\to} B(\gL')
  \]
  as $U_q'(\fg)$-crystals, and this isomorphism maps $u_{\gL} \otimes b$ to $u_{\gL'}$.
  {\normalfont(}Here both $B(\gL)$ and $B(\gL')$ are regarded as $U_q'(\fg)$-crystals. See Remark \ref{Rem:restriction}{\normalfont)}.
\end{Thm}

If $B$ is a perfect crystal of level $\ell$, then $\gee \circ \gph^{-1}$ 
induces a bijection $(P_\cl^+)^\ell \to (P_\cl^+)^\ell$, which is called the \textit{associated automorphism} of $B$.

For $i \in I$, we denote by $\tau^i \in \gS$
the unique element satisfying $t_{c_i \varpi_i} (\tau^i)^{-1} \in W$. 
Note that this definition is the same as that of Subsection \ref{subsection:Dynkin} for $i \in I^s$.
For $i \in I \setminus I^s$, $\tau^i$ are as follows: for $B_n^{(1)}$, $D_n^{(1)}$, $A_{2n-1}^{(2)}$, 
$\tau^i = \id$ if $i$ is even, and $\tau^i = \tau^1$ if $i$ is odd.
For $C_n^{(1)}$, $A_{2n}^{(2)}, D_{n+1}^{(2)}$, $\tau^i= \id$ for all $i \in I \setminus I^s$.

\begin{Thm}[\cite{MR2642564}]\label{Thm:perfectness} \ \\
  {\normalfont(i)} The level of a KR crystal $B^{r,\ell}$ is $\lceil \ell/c_r \rceil( = \min\{m \in \Z \mid m \ge 
                   \ell/c_r\})$, where $c_r$ is defined in Subsection \ref{subsection:Dynkin}. \\
  {\normalfont(ii)} $B^{r,\ell}$ is perfect if and only if $\ell/c_r \in \Z$. \\
  {\normalfont(iii)} The associated automorphism of $B^{r,c_r\ell}$ coincides with the action of $(\tau^r)^{-1}$ 
                     on $(P_\cl^+)^\ell$. 
\end{Thm}

\begin{proof}
  The assertions (i) and (ii) were proved in \cite{MR2642564}.
  The associated automorphism of each $B^{r,c_r\ell}$ is explicitly described in \cite{MR2642564},
  and we can check the assertion (iii) directly from them.
  We remark that the equation in \cite[Subsection 4.3]{MR2642564} for the associated automorphism $\tau$
  of $B^{n,\ell}$ for $D_n^{(1)}$ is misprint.
  It should be modified as follows:
  \[ \tau\big(\sum_{i = 0}^n \ell_i\gL_i\big) = \ell_{n}\gL_0 + \ell_{n-1} \gL_1 + \sum_{i = 2}^{n-2} \ell_i\gL_{n-i}
     + \begin{cases} \ell_0 \gL_{n-1} + \ell_1 \gL_n & n \ \text{even},\\
                     \ell_1 \gL_{n-1} + \ell_0 \gL_{n} & n \ \text{odd}.
       \end{cases}
  \]
\end{proof}

Let $B = B^{r,\ell}$ be a (not necessarily perfect) KR crystal.
$B$ is known to have a unique element belonging to $B_{\min}$, which we denote by $m(B)$, such that 
\[ \gee\big(m(B)\big) = \lev(B)\gL_0.
\]
(If $B$ is perfect, this fact follows from the definition. 
For non-perfect ones, see \cite[Lemma 3.11]{MR2945586}).
Similarly, $B$ has a unique element $m'(B) \in B_{\min}$ such that 
\[ \gph\big(m'(B)\big) = \lev(B)\gL_0.
\]
If $B$ is perfect, we have from Theorem \ref{Thm:perfectness} (iii) that
\begin{equation*} 
   \gph\big(m(B)\big) = \lev(B) \gL_{\tau^r(0)}.
\end{equation*}
The first assertion of the following theorem was proved in \cite{ST},
and the second one is obvious from Lemma \ref{Lem}:

\begin{Thm}[{\cite[Theorem 6.1]{ST}}] \label{Thm:Demazure}
  Let $B = B^{r,c_r\ell}$ be a perfect KR crystal.
  Then the isomorphism $B(\ell\gL_0) \otimes B \stackrel{\sim}{\to} B(\ell \gL_{\tau^r(0)})$
  given in Theorem \ref{Thm:Path} maps the subset $u_{\ell\gL_0} \otimes B$ onto the Demazure crystal 
  $ B_{t_{c_rw_0(\varpi_r)}}(\ell\gL_0)$.
  Moreover, the image of the element $u_{\ell\gL_0} \otimes u(B)$ under this isomorphism is the extremal 
  element $u_{t_{c_rw_0(\varpi_r)}(\ell\gL_0)}$.
\end{Thm}

Later we need the following lemma:

\begin{Lem} \label{Lem:elementary2}
  Let $B_1,B_2$ be perfect KR crystals, and assume that $\lev (B_1) \le \lev(B_2)$.
  If $b_1 \otimes b_2 \in (B_1 \otimes B_2)_{\min}$, then for every $b_2' \in B_2$ there exists a sequence
  $i_1,\dots,i_k$ of elements of $I$ such that
\begin{align} \label{eq:act_right}
  e_{i_k}\cdots e_{i_1}(b_1 \otimes b_2') &= b_1 \otimes (e_{i_k}\cdots e_{i_1} b_2') \nonumber \\
                                        &= b_1 \otimes b_2.
\end{align}
\end{Lem}

\begin{proof}
  By Lemma \ref{Lem:level} (iii), $b_2 \in (B_2)_{\min}$ and 
  \begin{equation}\label{eq:add}
    \gee(b_2) - \gph(b_1) \in P_{\cl}^+.
  \end{equation}
  Set $\gL = \aff\big(\gee(b_2)\big)$ and $\gL' = \aff\big( \gph(b_2)\big)$.
  By Theorem \ref{Thm:Path}, there exists an isomorphism
  \[ B(\gL) \otimes B_2 \stackrel{\sim}{\to}B(\gL')
  \]
  which maps $u_{\gL} \otimes b_2$ to $u_{\gL'}$.
  Therefore, there exists a sequence $i_1,\dots,i_k$ of elements of $I$ such that 
  \begin{align*}
    e_{i_k}\cdots e_{i_1}(u_{\gL} \otimes b_2') &= u_{\gL} \otimes (e_{i_k}\cdots e_{i_1} b_2') \\
                                                     &= u_{\gL} \otimes b_2.
  \end{align*}
  This equation implies for each $1 \le q \le k$ that 
  \[ \gee_{i_q}(e_{i_{q-1}} \cdots e_{i_1} b_2') > \gph_{i_q}(u_\gL) = \langle \cl(\gL), \ga_{i_q}^\vee\rangle 
     = \gee_{i_q}(b_2).
  \]
  Hence we have $\gee_{i_q}(e_{i_{q-1}} \cdots e_{i_1} b_2') >\gee_{i_q}(b_2) \ge \gph_{i_q}(b_1)$ 
  for each $1 \le q \le k$ by (\ref{eq:add}), and (\ref{eq:act_right}) is proved.
\end{proof}

\subsection{Isomorphism as full subgraphs of $U_q'(\fg)$-crystals}

We need the following elementary lemma:

\begin{Lem} \label{Lem:elementary}
  Let $B_1, B_2$ be crystals, and $b_j \in B_j$ $(j = 1, 2)$ arbitrary elements.
  If $f_i b_1 \neq 0$ for some $i \in I$, then there exist some $b_2'\in B_2$ and $m \in \Z_{>0}$ such that
  \[ f_ib_1 \otimes b_2 = f_i^m(b_1 \otimes b_2').
  \] 
\end{Lem}

\begin{proof}
  When $\gph_i(b_1) > \gee_i(b_2)$, $m = 1$ and $b_2' = b_2$ satisfy the assertion.
  When $\gph_i(b_1) \le \gee_i(b_2)$, $m = \gee_i(b_2) - \gph_i(b_1)+2$ and $b_2' = e_i^{m - 1} b_2$
  satisfy this.  
\end{proof}

Now we show the following proposition (the notion of an isomorphism of full subgraphs has been defined
in Section \ref{Section:crystals}):

\begin{Prop}\label{Prop:Main_Prop}
  Let $B_j = B^{r_j, c_{r_j}\ell_j}$ for $1 \le j \le p$ be perfect KR crystals with $\ell_1 \le \ell_2 \le \cdots \le \ell_p$,
  and set $\ell^j = \ell_j - \ell_{j-1}$ with $\ell_0 = 0$. 
  We put $\mu_j = c_{r_j}w_0(\varpi_{r_j})$ and $B = B_p \otimes \cdots \otimes B_2 \otimes B_1$.
  Then there exists an isomorphism 
  \begin{equation*} 
    \Psi_B\colon u_{\ell_p\gL_0} \otimes B \stackrel{\sim}{\to} 
    \cF_{t_{\mu_p}}\Big( u_{\ell^p \gL_0} \otimes \cdots \otimes \cF_{t_{\mu_2}}\Big(u_{\ell^2\gL_0} 
    \otimes \cF_{t_{\mu_1}}(u_{\ell^1\gL_0})\Big)\!\cdots\!\Big)
  \end{equation*}
  of full subgraphs of $U_q'(\fg)$-crystals, where the right hand side is a subset of
  $B\big(\ell^p\gL_{\tau^{p}(0)}\big) \otimes \cdots \otimes B\big(\ell^1\gL_{\tau^p\cdots \tau^1(0)}\big)$
  {\normalfont(}regarded as a $U_q'(\fg)$-crsytal\,{\normalfont)}.
  Moreover, the isomorphism $\Psi_B$ maps the element $u_{\ell_p\gL_0} \otimes u(B)$ to the tensor product of 
  the extremal elements 
  $u_{t_{\mu_p}(\ell^p\gL_0)} \otimes \cdots \otimes u_{t_{\mu_p + \cdots + \mu_2}(\ell^2\gL_0)} \otimes 
     u_{t_{\mu_p + \cdots + \mu_1} (\ell^1\gL_0)}.$ 
\end{Prop}

\begin{proof}
  If $p = 1$, the assertion follows from Theorem \ref{Thm:Demazure}.
  Assume $p > 1$. We put $\tau = \tau^{r_p}$ and $w = t_{\mu_p}\tau^{-1} \in W$. 
  Since $u_{\ell_p \gL_0} \otimes B_p \stackrel{\sim}{\to} \cF_{w}(u_{\ell_p\gL_{\tau(0)}})$, 
  we have
  \[ u_{\ell_p\gL_0} \otimes B_p \otimes \cdots \otimes B_1 \stackrel{\sim}{\to}
     \cF_{w}(u_{\ell_p\gL_{\tau(0)}}) \otimes B_{p-1} \otimes \cdots \otimes B_1,
  \]
  and the right hand side is equal to 
  $\cF_{w}(u_{\ell_p\gL_{\tau(0)}}\otimes B_{p-1} \otimes \cdots\otimes B_1)$ by Lemma \ref{Lem:elementary}.
  Then since it follows from Proposition \ref{Prop:ST} that
  \begin{align*}
    u_{\ell_p \gL_{\tau(0)}} \otimes B_{p-1} \otimes \cdots \otimes &B_1  \cong \wti{\tau}\left(u_{\ell_p \gL_0} 
    \otimes B_{p-1} \otimes \cdots \otimes B_1\right) \\
    & \cong \wti{\tau}\Big(u_{\ell^p\gL_0} \otimes \big(u_{\ell_{p-1} \gL_0} \otimes B_{p-1} 
      \otimes \cdots \otimes B_1\big)\Big),
  \end{align*}
  we obtain an isomorphism $\Psi_B$ by the induction hypothesis.  
  For each $1 \le j \le p$, we have  $t_{\mu_p} \cdots t_{\mu_j} = t_{\mu_p + \cdots + \mu_j}$
  and $\ell(t_{\mu_p+ \cdots + \mu_j}) = \sum_{k=j}^p \ell(t_{\mu_k})$
  since $\ell(w) = \#\{ \ga \in \gD \cap Q^+ \mid w(\ga) \in - Q^+\}$ for $w \in \wti{W}$.
  Therefore it follows from Proposition \ref{Prop:twisted_Demazure} that the right hand side of $\Psi_B$ is 
  contained in the set
  \[ B_{t_{\mu_p}}(\ell^p\gL_0) \otimes\cdots \otimes B_{t_{\mu_p + \cdots +\mu_2}}(\ell^2\gL_0) 
     \otimes B_{t_{\mu_p + \cdots + \mu_1}}(\ell^1\gL_0)
  \]
  by definition.
  Then the second assertion is easily proved from Lemma \ref{Lem}, since we have
  \begin{align*}
     \langle t_{\mu_p + \cdots +\mu_j}( \ell^j \gL_0), \ga_i^\vee \rangle 
     &= \ell^j\sum_{k=j}^p c_{r_k}\big\langle w_0(\varpi_{r_k}), \ga_i^{\vee} \big\rangle \le 0 \ \text{for all} \
      i \in I_0,  \ \ \ \text{and}    \\
     \cl \Big(\sum_{j=1}^p t_{\mu_p + \cdots + \mu_j}(\ell^j \gL_0) \Big)
     &= \sum_{j=1}^p c_{r_j}\ell_j w_0(\varpi_{r_j}) = \wt\big(u(B)\big).
  \end{align*}
\end{proof}

\begin{Rem} \label{Rem}\normalfont
(i)  Put $B^{p-1} = B_{p-1} \otimes \cdots \otimes B_2 \otimes B_1$.
     We see from the construction of the isomorphism $\Psi_B$ that the following diagram commutes 
     (where we set $\tau = \tau^{r_p}$):
  \[ \xymatrix{ u_{\ell_{p-1}\gL_0} \otimes B_{p-1} \otimes \cdots \otimes B_1 \ar[r]_{\Psi_{B^{p-1}}}^(.48){\cong} 
     \ar[d]^{\varphi} &  
     \cF_{t_{\mu_{p-1}}}(u_{\ell^{p-1}\gL_0} \otimes \cdots)  \ar[d]^{\psi}  \\
     u_{\ell_p\gL_0} \otimes m(B_p) \otimes B_{p-1} \otimes \dots \otimes B_1 \ar@{^{(}->}[d] & 
     \widetilde{\tau}\big(u_{\ell^p\gL_0} \otimes \cF_{t_{\mu_{p-1}}}(u_{\ell^{p-1}\gL_0} \otimes\cdots)\big) \ar@{^{(}->}[d] \\ 
     u_{\ell_p\gL_{0}} \otimes B_p \otimes B_{p-1} \otimes\cdots \otimes B_1 \ar[r]_(.47){\Psi_B}^(.46){\cong} &
     \cF_{t_{\mu_p}}\big(u_{\ell^p\gL_0} \otimes \cF_{t_{\mu_{p-1}}}(u_{\ell^{p-1}\gL_0} \otimes\cdots)\big),
     }
  \]  
  where $\varphi$ and $\psi$ are bijective maps of sets defined by
  \[ \varphi(u_{\ell_{p-1}\gL_0} \otimes b) =u_{\ell_p\gL_0} \otimes m(B_p) \otimes \tau(b) \ \ \ 
     \text{for} \ b \in B^{p-1}, \ \ \ \text{and}
  \]
  \[ \psi(b) = \wti{\tau}\big(u_{\ell^p\gL_0} \otimes b\big) \ \ \ \text{for} \ 
     b \in \cF_{t_{\mu_{p-1}}}(u_{\ell^{p-1}\gL_0} \otimes \cdots)
  \]
  respectively. \\
(ii) By Proposition \ref{Prop:twisted_Demazure} and Theorem \ref{Thm:excellent_filtration},
    the right hand side of the isomorphism $\Psi_B$ is isomorphic as a full subgraph to a disjoint union of Demazure crystals.\\
(iii) The right hand side of $\Psi_B$ also appeared in \cite{MR1887117}
    as the crystal basis of a \textit{generalized Demazure module}.
\end{Rem}

Note that the right hand side of $\Psi_B$ is a subset of a tensor product of the crystal bases of 
$U_q(\fg)$-modules.
Hence each element $b$ of this set has a natural $\Z$-grading given by $\langle \wt(b), d \rangle$.
The goal of this article is to show that, 
under the isomorphism $\Psi_B$, the minus of this natural grading coincides up to a shift 
with the grading on the left hand side given by the energy function introduced in the next section.

\section{Energy function}

Similarly as \cite{MR1187560}, the following proposition is proved from the existence of the universal $R$-matrix
and Theorem \ref{Thm:KR}:

\begin{Prop} Let $B_1, B_2 \in \cC$. \\
  {\normalfont(i)} There exists a unique isomorphism 
    $\gs=\gs_{B_1,B_2}\colon B_1 \otimes B_2 \stackrel{\sim}{\to} B_2 \otimes B_1$ 
    of $U_q'(\fg)$-crystals called the \textit{combinatorial $R$-matrix}. \\
  {\normalfont(ii)} There exists a unique map $H = H_{B_1,B_2}\colon B_1 \otimes B_2 \to \Z$ 
    called the \textit{local energy function}
    such that $H\big(u(B_1\otimes B_2)\big) = 0$, $H$ is constant on each $U_q(\fg_0)$-component, and for $b_1 \otimes b_2
    \in B_1 \otimes B_2$ mapped to $\widetilde{b}_2 \otimes \widetilde{b}_1 \in B_2 \otimes B_1$ under $\gs$, we have
    \begin{align*}
      H\big(&e_0(b_1 \otimes b_2)\big) \\ 
      &= \begin{cases} H(b_1 \otimes b_2) + 1 & \text{if} \ e_0(b_1 \otimes b_2) = e_0b_1 \otimes b_2, \
                     e_0(\widetilde{b}_2 \otimes \widetilde{b}_1) = e_0 \widetilde{b}_2 \otimes \widetilde{b}_1, \\
                     H(b_1 \otimes b_2) - 1 & \text{if} \ e_0(b_1 \otimes b_2) = b_1 \otimes e_0b_2, \ 
                     e_0(\widetilde{b}_2 \otimes \widetilde{b}_1) = \widetilde{b}_2 \otimes e_0\widetilde{b}_1, \\
                     H(b_1 \otimes b_2) & \text{otherwise}.
         \end{cases}
    \end{align*}                                                                        
\end{Prop}

For $B_1, B_2 \in \cC$, we have $\gs\big(u(B_1) \otimes u(B_2)\big) = u(B_2) \otimes u(B_1)$ 
by the weight consideration.
Recall that for every $\tau \in \gS$, there exists some $w\in W_0$ such that $\tau\big(u(B_1) \otimes u(B_2)\big)
=S_w\big(u(B_1) \otimes u(B_2)\big)$ by Lemma \ref{Lem:highest}.
Hence we have
\begin{align*}
   \gs \circ \tau\big(u(B_1) \otimes u(B_2)\big) &= \gs \circ S_w \big(u(B_1) \otimes u(B_2)\big) \\ 
   = S_w\big(u(B_2) \otimes u(B_1)\big) &= \tau\big(u(B_2) \otimes u(B_1)\big),
\end{align*}
which together with (\ref{eq:tau}) implies that $\gs$ commutes with the action of $\tau$.
The following lemma is a consequence of the definition of the local energy function:

\begin{Lem}\label{Lem:local_energy_function}
  Let $B_1,B_2 \in \cC$, $b_j \in B_j$ for $j=1,2$ be such that $\gs(b_1 \otimes b_2) = \wti{b}_2 \otimes \wti{b}_1$,
  and $j_1,\ldots, j_\ell$ an arbitrary sequence of elements of $I$ satisfying 
  $e_{j_\ell} \cdots e_{j_1} (b_1 \otimes b_2) \neq 0$.
  If
  \begin{align*} 
    e_{j_\ell} \cdots e_{j_1} (b_1 \otimes b_2) &= e_{i_{\ell-k}'} \cdots e_{i_1'} b_1 \otimes e_{i_k} \cdots e_{i_1}b_2 \ \
    \text{and} \\
    e_{j_\ell} \cdots e_{j_1} (\wti{b}_2 \otimes \wti{b}_1) 
     &= e_{\tilde{i}_{m}} \cdots e_{\tilde{i}_1}\wti{b}_2 \otimes e_{\tilde{i}'_{\ell-m}} \cdots e_{\tilde{i}'_1} \wti{b}_1
  \end{align*}
  hold where 
  \[ \{j_1, \ldots, j_\ell \} = \{i_1, \ldots, i_k\} \sqcup \{i_1',\ldots,i_{\ell-k}'\} = \{\tilde{i}_1,\ldots,\tilde{i}_m\} 
     \sqcup \{ \tilde{i}_1', \ldots, \tilde{i}_{\ell-m}'\}
  \]
  as multisets,
  then we have
  \begin{align*} H\big(e_{j_\ell} \cdots e_{j_1} (b_1 \otimes b_2)\big) &- H(b_1 \otimes b_2) \\
    &=\#\{1 \le q \le m \mid \tilde{i}_q = 0\} - \#\{1 \le q \le k \mid i_q = 0 \}.
  \end{align*}
\end{Lem}

For $B \in \cC$, the \textit{energy function} $D = D_B\colon B \to \Z$ is defined as follows: \\
(i) If $B$ is a single KR crystal, then define
    \[ D_B(b) = H_{B,B}\big(m'(B) \otimes b\big) -H_{B,B}\big(m'(B) \otimes u(B)\big).
    \]
(ii) If $B_1,B_2 \in \cC$ and $B = B_1 \otimes B_2$, then define
     \[ D_B(b_1 \otimes b_2) = D_{B_1}(b_1) + D_{B_2}(\widetilde{b}_2) + H_{B_1,B_2}(b_1 \otimes b_2),
     \]
     where $\gs_{B_1, B_2}(b_1 \otimes b_2) = \widetilde{b}_2 \otimes \widetilde{b}_1$.\\
     
By definition, $D_B$ is constant on each $U_q(\fg_0)$-component of $B$ and we have
\begin{equation}\label{eq:0}
  D_B\big(u(B)\big) = 0.
\end{equation}

\begin{Prop}[\cite{MR1973369}] \label{Prop:energy_function}
  {\normalfont(i)} For $B_1, B_2, B_3 \in \cC$, we have
       \[ D_{(B_1 \otimes B_2)\otimes B_3} = D_{B_1 \otimes (B_2 \otimes B_3)}.
       \]
       Hence for every $B \in \cC$, the function $D_B$ is well-defined. \\
  {\normalfont(ii)} Let $B = B_1 \otimes \cdots \otimes B_p \in \cC$.
       For $b_1 \otimes \cdots \otimes b_p \in B$ and $1 \le i \le j \le p$, 
       define $b_j^{(i)} \in B_j$ by 
       \begin{align*}   
         B_i \otimes \cdots \otimes B_{j-1} \otimes B_j &\stackrel{\sim}{\to} B_j \otimes B_i \otimes \cdots \otimes B_{j-1} \\
         b_i \otimes \cdots \otimes b_{j-1} \otimes b_j &\mapsto b_j^{(i)} \otimes \wti{b}_i 
         \otimes \cdots \otimes \wti{b}_{j-1}.
       \end{align*}
       Then we have 
       \[ D_B(b_1\otimes \cdots \otimes b_p) = \sum_{1 \le j \le p} D_{B_j}(b_j^{(1)}) + \sum_{1 \le j < k \le p} 
          H_{B_j,B_k}(b_j \otimes b_k^{(j+1)}).
       \]
\end{Prop}

\begin{Lem}\label{Lem:elementary_for_energy}
  Let $B \in \cC$ and $\ell = \lev (B)$.
  If $b \in B$ satisfies $\gee_0(b) > \ell$, then we have
  \[ D(e_0b) = D(b) -1.
  \]
\end{Lem}

\begin{proof}
  We show the assertion by the induction on the number $p$ of tensor factors of $B$.
  The case $p=1$ follows since we have
  \[ D(e_0b) = H\big(m'(B) \otimes e_0b\big) - t = H\big(m'(B) \otimes b \big) -1
     -t = D(b) -1,
  \]
  where we set $t =H\big(m'(B) \otimes u(B)\big)$.
  Assume $p > 1$, and write $B = B_1 \otimes B_2$, $b = b_1 \otimes b_2$ and 
  $\wti{b}_2 \otimes \wti{b}_1 = \gs(b_1 \otimes b_2)$. 
  We can show the assertion by computing case by case, using $\lev(B_1) \le \ell$ and $\lev(B_2) \le \ell$. 
  For example, assume that $e_0(b_1 \otimes b_2) = e_0b_1 \otimes b_2$ and 
  $e_0(\wti{b}_2 \otimes \wti{b}_1) = e_0\wti{b}_2 \otimes \wti{b}_1$.
  Then we have $\gee_0(b_1) = \gee_0(b_1 \otimes b_2) > \ell$
  and $\gee_0(\wti{b}_2) = \gee_0(\wti{b}_2 \otimes \wti{b}_1) >\ell$,
  which imply by the induction hypothesis that
  \begin{align*}
    D(e_0b) &= D(e_0b_1) + D(e_0\wti{b}_2) + H\big(e_0(b_1 \otimes b_2)\big) \\
            &= \big(D(b_1) -1\big) + \big(D(\wti{b}_2) -1 \big)+ \big(H(b_1 \otimes b_2) + 1\big)\\
            &= D(b) -1.
  \end{align*}
  The other cases are proved similarly.
 \end{proof}

\section{Main theorem}

\subsection{Statement and corollaries}

Now, we state the main theorem of this article.
This theorem is a generalization of \cite[Theorem 7.4]{ST}, in which $\ell_1 = \ell_2 = \cdots =\ell_p$ is assumed.

\begin{Thm}\label{Thm:Main_Theorem}
  Let $B_j = B^{r_j, c_{r_j}\ell_j}$ for $1 \le j \le p$ be perfect KR crystals with $\ell_1 \le \ell_2 \le \cdots \le \ell_p$,
  and set $\ell^j = \ell_j - \ell_{j-1}$ with $\ell_0 = 0$. 
  We put $\mu_j = c_{r_j}w_0(\varpi_{r_j})$ and $B = B_p \otimes \cdots \otimes B_2 \otimes B_1$.
  Then there exists an isomorphism 
  \begin{equation*} 
    \Psi_B\colon u_{\ell_p\gL_0} \otimes B \stackrel{\sim}{\to} 
    \cF_{t_{\mu_p}}\Big( u_{\ell^p \gL_0} \otimes \cdots \otimes \cF_{t_{\mu_2}}\Big(u_{\ell^2\gL_0} 
    \otimes \cF_{t_{\mu_1}}(u_{\ell^1\gL_0})\Big)\!\cdots\!\Big)
  \end{equation*}
  of full subgraphs of $U_q'(\fg)$-crystals satisfying
  \begin{equation}\label{eq:weight_id}
    D(b) = - \langle \wt\, \Psi_B(u_{\ell_p \gL_0} \otimes b), d \rangle + C_B
  \end{equation}
  for every $b \in B$, where $C_B \in N^{-1}\Z$ is a constant defined by
  \[ C_B= \sum_{j=1}^p \ell^j \big\langle t_{\mu_p + \cdots + \mu_{j+1} + \mu_j}(\gL_0), d \big\rangle.
  \]
\end{Thm}

Recall that, as stated in Remark \ref{Rem} (ii), the right hand side of $\Psi_B$ is isomorphic 
as a full subgraph to a disjoint union of Demazure crystals.
Hence we can see inductively using Corollary \ref{Cor:Demazure_character} that the following equation holds:

\begin{Cor} \label{Cor:Main_cor}
  Under the notation and the assumptions of Theorem \ref{Thm:Main_Theorem}, we have
  \begin{align*}
    e^{\ell_p\gL_0+ C_B\gd} \sum_{b \in B} e^{\aff\circ\wt(b) - \gd D(b)}
    = D_{t_{\mu_p}}\left( e^{\ell^p\gL_0} \cdots D_{t_{\mu_2}}
    \left( e^{\ell^2\gL_0} \cdot D_{t_{\mu_1}}(e^{\ell^1\gL_0})\right)\!\cdots\! \right).
  \end{align*}
\end{Cor}

As \cite{MR1745263,MR1903978}, the \textit{one-dimensional sum} $X(B, \mu, q) \in \Z[q, q^{-1}]$
for $\mu \in P_0^+$ is defined by
\[ X(B, \mu ,q) = \sum_{\begin{smallmatrix} b \in B \\ e_i b= 0 \ (i \in I_0) \\ \wt(b) = \mu \end{smallmatrix}} q^{D(b)}.
\]
Let $\ch V_{\fg_0}(\mu)$ denote the character of the irreducible $\fg_0$-module with highest weight $\mu$.
Since 
\[ \sum_{b \in B} q^{D(b)} e^{\wt(b)} = \sum_{\mu \in P_0^+} X(B,\mu,q) \ch V_{\fg_0}(\mu)
\]
holds, we have the following corollary:

\begin{Cor}\label{Cor:X}
  Under the notation and the assumptions of Theorem \ref{Thm:Main_Theorem}, we have
  \begin{align*}
    q^{-C_B}\sum_{\mu \in P_0^+}& X(B, \mu,q) \ch V_{\fg_0}(\mu) \\ 
    &= e^{-\ell_p\gL_0}D_{t_{\mu_p}}\Big( e^{\ell^p\gL_0} \cdots D_{t_{\mu_2}}
    \left( e^{\ell^2\gL_0} \cdot D_{t_{\mu_1}}(e^{\ell^1\gL_0})\Big)\!\cdots\! \right),
  \end{align*}
  where we set $q = e^{-\gd}$ and consider $\ch V_{\fg_0}(\mu)$ as an element of $\C[P]$ via the map $\aff\colon P_\cl \to P$.
\end{Cor}

\begin{Rem} \normalfont
  Let $\eta $ be a permutation of the set $\{1, \ldots, p\}$, 
  and put $B_\eta= B_{\eta(1)} \otimes \cdots \otimes B_{\eta(p)}$.
  We have from \cite[Lemma 2.15]{MR1973369} that
  \[ D_{B_{\eta}}\big(\gs_{\eta}(b)\big) = D_B(b) \ \ \ \text{for every} \ b \in B,
  \]
  where $\gs_{\eta}\colon B \stackrel{\sim}{\to} B_{\eta}$ is the unique isomorphism.
  In particular, we have 
  \[ X(B,\mu,q) = X(B_{\eta}, \mu, q).
  \]
  Hence for every $\eta$, the above theorem and corollaries with $B$ replaced by $B_{\eta}$ 
  (and the right hand sides unchanged) also hold. 
\end{Rem}

\subsection{Proof of the main theorem}

In order to prove the main theorem, it remains to show that the isomorphism $\Psi_B$ constructed 
in Proposition \ref{Prop:Main_Prop} satisfies 
(\ref{eq:weight_id}).
To show this, we prepare several lemmas.

\begin{Lem}\label{Lem:Lemma1}
  Let $B_1, B_2 \in \cC$ and $\tau \in \gS$.
  For $b_1 \otimes b_2 \in B_1 \otimes B_2$ mapped to $\tilde{b}_2 \otimes \tilde{b}_1 \in B_2 \otimes B_1$ 
  under $\gs$, we have
  \begin{equation}\label{eq:difference}
    H\big(b_1 \otimes b_2\big) - H\big(\tau(b_1 \otimes b_2)\big)
    = \langle \wt (b_2) - \wt (\tilde{b}_2), \varpi_{\tau^{-1}(0)}^{\vee} \rangle.
  \end{equation}
\end{Lem}

\begin{proof}
  Although the proof is carried out in a similar way as that of \cite[Lemma 8.2]{MR2945586},
  we give it for the reader's convenience.
  
  The case $\tau = \id$ is trivial. We assume otherwise, and put $t = \tau^{-1}(0) \in I^s \setminus \{ 0\}$.
  If $b_1 = u(B_1)$ and $b_2 = u(B_2)$, we have from Lemma \ref{Lem:highest} that
  \[ H\big(\tau(u(B_1) \otimes u(B_2))\big) = H\big(u(B_1) \otimes u(B_2)\big) = 0,
  \]
  and hence the left hand side of (\ref{eq:difference}) is $0$.
  On the other hand, the right hand side is also $0$
  since we have $\gs\big(u(B_1) \otimes u(B_2)\big) = u(B_2) \otimes u(B_1)$, and the assertion is proved in this case.
  Therefore by Lemma \ref{Lem:strongly_connected}, it suffices to show for each $i \in I$ 
  that if (\ref{eq:difference}) holds and $e_i(b_1 \otimes b_2) \neq 0$,
  then (\ref{eq:difference}) with $b_1 \otimes b_2$ replaced by $e_i(b_1 \otimes b_2)$ also holds. 
  If $i \neq 0, t$, it is easy to see that the both sides of (\ref{eq:difference}) 
  do not change when $b_1 \otimes b_2$ is replaced by $e_i(b_1 \otimes b_2)$.
  Assume that $i = 0$.
  Since $t \neq 0$, we have
  \begin{align} \label{align}
    &\Big( H\big(e_0(b_1 \otimes b_2)\big) - H\big(\tau \circ e_0(b_1 \otimes b_2)\big) \Big) - 
    \Big( H\big(b_1 \otimes b_2\big) - H\big(\tau(b_1 \otimes b_2)\big)\Big) \nonumber \\
    &= \Big(H\big( e_0(b_1 \otimes b_2)\big) - H(b_1 \otimes b_2)\Big) - \Big(H\big(e_t\circ \tau(b_1\otimes b_2)\big)
            - H\big(\tau(b_1\otimes b_2)\big)\Big) \nonumber \\
    &= \begin{cases} 1 & \text{if} \ e_0(b_1 \otimes b_2) = e_0b_1 \otimes b_2, \
                     e_0(\widetilde{b}_2 \otimes \widetilde{b}_1) = e_0 \widetilde{b}_2 \otimes \widetilde{b}_1, \\
                      - 1 & \text{if} \ e_0(b_1 \otimes b_2) = b_1 \otimes e_0b_2, \ 
                     e_0(\widetilde{b}_2 \otimes \widetilde{b}_1) = \widetilde{b}_2 \otimes e_0\widetilde{b}_1, \\
                     0 & \text{otherwise}.
         \end{cases}
  \end{align}
  On the other hand, putting $e_0(b_1 \otimes b_2) = b_1' \otimes b_2'$ and $e_0(\tilde{b}_2 \otimes \tilde{b}_1) 
  = \tilde{b}_2' \otimes \tilde{b}_1'$, we easily check using $\langle \ga_0, \varpi_t^\vee\rangle = -1$ (see (\ref{eq:pairing})) that 
  $\langle \wt (b_2') - \wt (\tilde{b}_2'), \varpi_{t}^{\vee} \rangle - 
     \langle \wt(b_2) - \wt (\tilde{b}_2), \varpi_{t}^{\vee} \rangle$
  is equal to (\ref{align}), 
  which implies the assertion for $i = 0$. The case $i = t$ is similar.
\end{proof}

For $B \in \cC$ and $\ell \in \Z_{>0}$, we define a subset $\hw_{I_0}^{\le \ell}(B) \subseteq B$ by
\[ \hw_{I_0}^{\le \ell}(B) = \{ b \in B \mid b \ \text{is $U_q(\fg_0)$-highest weight, } \gee_0(b) \le \ell \}.
\]

\begin{Lem}\label{Lem:Lemma2}
  Let $B_j = B^{r_j,c_{r_j}\ell_j} (j = 1,2)$ be two perfect KR crystals, and assume that $\ell_1 \ge \ell_2$.
  For every $b_2 \in \hw_{I_0}^{\le \ell_1}(B_2)$, we have
  \[ \gs_{B_1, B_2}\big(m(B_1) \otimes \tau^{r_1}(b_2) \big)= b_2 \otimes b_1
  \]
  for some $b_1 \in B_1$.
\end{Lem}
  
\begin{proof}
  Put $\tau = \tau^{r_1}$.
  Since 
  \[ \gph\big(m(B_1)\big) = \ell_1\gL_{\tau(0)} \  \text{and} \ \gee\big(\tau(b_2)\big) = \gee_0(b_2) \gL_{\tau(0)},
  \]
  we have from Lemma \ref{Lem:level} (ii) that $m(B_1) \otimes \tau(b_2) \in (B_1 \otimes B_2)_{\min}$, 
  $\gee\big(m(B_1) \otimes \tau(b_2)\big) = \ell_1 \gL_0$, and
  \[ \gph\big(m(B_1) \otimes \tau(b_2)\big) = \big(\ell_1 - \gee_0(b_2)\big)\gL_{\tau(0)} + \tau\big(\gph(b_2)\big).
  \] 
  Put $\gL = \big(\ell_1 - \gee_0(b_2)\big)\gL_{\tau(0)} + \tau\big(\gph(b_2)\big)$, 
  and $b_2' \otimes b_1' = \gs\big(m(B_1) \otimes \tau(b_2)\big)$.
  Since $b_2' \otimes b_1' \in (B_2 \otimes B_1)_{\min}$,
  we have from Lemma \ref{Lem:level} (iii) that $b_1' \in (B_1)_{\min}$ and 
  $\gph(b_1') = \gph(b_2' \otimes b_1') = \gL$.
  Hence from Theorem \ref{Thm:perfectness} (iii), we have 
  \begin{equation}\label{eq:gee}
    \gee(b_1') = \tau^{-1}(\gL) =\big(\ell_1 -\gee_0(b_2)\big) \gL_0 +\gph(b_2).
  \end{equation}
  Note that $\gee(b_2' \otimes b_1') = \gee\big(m(B_1) \otimes \tau(b_2)\big) = \ell_1 \gL_0$,
  and this also implies $\gee(b_2') = \ell\gL_0$ with $\ell \le \ell_1$ by (\ref{eq:tensor}).
  Hence we have that
  \begin{align*}
    \gph(b_2') = \wt(b_2') + \gee(b_2') = \gee(b_1') -\gee(b_2' \otimes b_1') + \gee(b_2') \in \gph(b_2) + \Z \gL_0,
  \end{align*}
  where the second equality follows from Lemma \ref{Lem:level} (iii).
  Then we have $b_2 = b_2'$ by Corollary \ref{Cor:multiplicity_free}, as required.
\end{proof}

\begin{Lem}\label{Lem:Lemma3}
  Let $B_j = B^{r_j,c_{r_j}\ell_j} (j = 1,2)$ be two perfect KR crystals, and assume that $\ell_1 \ge \ell_2$.
  Then there exists some global constant $C$ such that 
  \[ H\big(m(B_1) \otimes \tau(b_2)\big) = -\langle \wt(b_2), \varpi_{\tau^{-1}(0)}^{\vee} \rangle + C
  \]
  for every $b_2 \in \hw_{I_0}^{\le \ell_1}(B_2)$, where we put $\tau = \tau^{r_1}$.
\end{Lem}

\begin{proof}
  Although the proof of this lemma is basically the same as that of \cite[Lemma 4.7]{MR1993475},
  we include it for the reader's convenience.
  
  It suffices to show for $b_2,b_2^\dag \in \hw_{I_0}^{\le \ell_1}(B_2)$ that
  \[ H\big(m(B_1) \otimes \tau(b_2^\dag)\big) - H\big(m(B_1) \otimes \tau(b_2)\big) = - \langle \wt (b_2^\dag) - \wt (b_2), 
     \varpi_{\tau^{-1}(0)}^{\vee} \rangle.
  \]
  By Lemma \ref{Lem:Lemma2}, we have
  \[ \gs\big(m(B_1) \otimes \tau(b_2) \big) = b_2 \otimes b_1 \ \text{and} \ \gs\big(m(B_1) \otimes \tau(b_2^\dag) \big)
     = b_2^{\dag} \otimes {b_1^{\dag}}
  \]
  for some $b_1, {b_1^{\dag}} \in B_1$.
  By Lemma \ref{Lem:strongly_connected}, there exists a sequence $i_1,\ldots,i_k$ of elements of $I$ such that
  \[ e_{i_k} \cdots e_{i_1} b_2 = b_2^\dag.
  \]
  We choose such a sequence so that $k$ is minimal. 
  Using Lemma \ref{Lem:add}, we can take a sequence $j_1,\ldots,j_\ell \in I$ satisfying
  \begin{align*}
    e_{j_\ell} \cdots e_{j_1}(b_2 \otimes b_1) &= e_{i_k} \cdots e_{i_1}b_2 
                                           \otimes e_{i_{\ell - k}'} \cdots e_{i_1'}b_1 \\
                                        &= b_2^{\dag} \otimes e_{i_{\ell - k}'} \cdots e_{i_1'}b_1.
  \end{align*}
  Since $b_2^{\dag} \otimes {b_1^{\dag}} \in (B_2 \otimes B_1)_{\min}$, we may assume
  $e_{i_{\ell - k}'} \cdots e_{i_1'}b_1 = {b_1^{\dag}}$ by Lemma \ref{Lem:elementary2}.
  Then we have 
  \[ e_{j_\ell}\cdots e_{j_1} \big(m(B_1) \otimes \tau(b_2)\big) = m(B_1) \otimes \tau(b_2^{\dag}).
  \]
  We define the two sequences $\tilde{i}_1, \ldots, \tilde{i}_m$ and $\tilde{i}_1', \ldots, \tilde{i}_{\ell-m}'$ of elements 
  of $I$ by
  \begin{align*}
    e_{j_\ell}\cdots e_{j_1} \big(m(B_1) \otimes \tau(b_2)\big) &= e_{\tilde{i}_{\ell -m}'} \cdots e_{\tilde{i}_1'}m(B_1)
    \otimes e_{\tilde{i}_m} \cdots e_{\tilde{i}_1}\tau(b_2)\\
    &= m(B_1) \otimes \tau(b_2^{\dag}).
  \end{align*}
  Since 
  \[ e_{\tilde{i}_m} \cdots e_{\tilde{i}_1}\tau(b_2) = \tau(e_{\tau^{-1}(\tilde{i}_m)} \cdots 
     e_{\tau^{-1}(\tilde{i}_1)}b_2) = \tau(b_2^{\dag}),
  \]
  we have $e_{\tau^{-1}(\tilde{i}_m)} \cdots 
     e_{\tau^{-1}(\tilde{i}_1)}b_2 = b_2^{\dag}$, which implies 
  \begin{equation}\label{eq:inequality1}
    \sum_{1 \le q \le m} \ga_{\tau^{-1}(\tilde{i}_q)} - \sum_{1 \le q \le k} \ga_{i_q} \in \Z_{\ge 0} \gd
  \end{equation}
  by the minimality of $k$.
  
  By repeating the above procedure interchanging the roles of $b_2$ and $b_2^{\dag}$,
  we obtain sequences of elements of $I$ satisfying the following:
  \begin{align}\label{eq:inequality2}
    e_{j_{\ell^*}^*} \cdots e_{j_1^*} (b_2^{\dag} \otimes b_1^\dag) &=e_{i_{k^*}^*} \cdots e_{i_1^*}b_2^\dag
                       \otimes e_{{i'}^*_{\ell^* - k^*}} \cdots e_{{i'_1}^*}b_1^{\dag} \nonumber \\
                       &= b_2 \otimes b_1, \nonumber \\
     e_{j_{\ell^*}^*} \cdots e_{j_1^*} \big(m(B_1) \otimes \tau(b_2^{\dag})\big) &=e_{{\tilde{i'}}^*_{\ell^* - m^*}} 
                       \cdots e_{\tilde{i'_1}^*}m(B_1) \otimes e_{\tilde{i}_{m^*}^*} \cdots 
                       e_{\tilde{i}_1^*}\tau(b_2^{\dag}) \nonumber \\
                       &= m(B_1) \otimes \tau(b_2), \nonumber \\
     \sum_{1 \le q \le m^*} \ga_{\tau^{-1}(\tilde{i}_q^*)} &- \sum_{1 \le q \le k^*} \ga_{i_q^*} \in \Z_{\ge 0} \gd.
  \end{align}
  By Lemma \ref{Lem:local_energy_function}, we have
  \begin{align*}
    0 &= \Big(H\big(m(B_1) \otimes \tau(b_2^{\dag})\big) - H\big(m(B_1) \otimes \tau(b_2)\big)\Big) \\
      & \hspace{40pt}+ \Big( H\big(m(B_1) \otimes \tau(b_2)\big) - H\big(m(B_1) \otimes \tau(b_2^{\dag})\big)\Big) \\
      &=\Big(\#\{1\le q \le k\mid i_q = 0\} - \#\{1 \le q \le m \mid \tilde{i}_q =0\}\Big) \\
      & \hspace{40pt} + \Big(\#\{1\le q \le k^*\mid i^*_q = 0\} - \#\{1 \le q \le m^* \mid \tilde{i}^*_q =0\}\Big)\\
      &=\Big(\#\{1\le q \le k\mid i_q = 0\} + \#\{1\le q \le k^*\mid i^*_q = 0\} \Big) \\
      & \hspace{40pt} - \Big(\#\{1 \le q \le m \mid \tilde{i}_q =0\} + \#\{1 \le q \le m^* \mid \tilde{i}^*_q =0\}\Big).
  \end{align*}
  This equation together with 
  (\ref{eq:inequality1}) and (\ref{eq:inequality2}) implies that $k =m$, $k^* = m^*$, and
  \begin{align*}
    \{\tau^{-1}(\tilde{i}_1), \ldots, \tau^{-1}(\tilde{i}_m)\} = \{i_1,\ldots,i_k\}, \ \
    \{\tau^{-1}(\tilde{i}_1^*), \ldots, \tau^{-1}(\tilde{i}^*_{m^*})\} = \{ i^*_1,\ldots,i_{k^*}^* \}
  \end{align*}
  as multisets.
  Hence we have
  \begin{align*}
    &H\big(m(B_1) \otimes \tau(b_2^\dag)\big) - H\big(m(B_1) \otimes \tau(b_2)\big)\\
    =& \#\{1\le q \le k\mid i_q = 0\} - \#\{1 \le q \le m \mid \tilde{i}_q =0\} \\
    =& \#\{1\le q \le k\mid i_q = 0\} - \#\{ 1 \le q \le k \mid i_q = \tau^{-1}(0)\} \\
    =& - \langle \wt(b_2^\dag) - \wt(b_2), \varpi_{\tau^{-1}(0)}^{\vee} \rangle,
  \end{align*}
  and the assertion is proved.
\end{proof}

The following lemma is crucial for the proof of our theorem:

\begin{Lem}\label{Lem:important}
  Let $B_j = B^{r_j,c_{r_j}\ell_j} (0 \le j \le p)$ be perfect KR crystals, and put
  $B = B_1 \otimes B_2 \otimes \cdots \otimes B_p$.
  We assume that $\ell_0 \ge \ell_j$ for every $1 \le j \le p$.
  Then there exists some global constant $C$ such that
  \[ D\big(m(B_0) \otimes \tau(b)\big) = D(b) - \big\langle \wt(b), \varpi_{\tau^{-1}(0)}^\vee \big\rangle + C
  \]
  for every $b \in \hw_{I_0}^{\le \ell_0}(B)$, where we put $\tau = \tau^{r_0}$.
\end{Lem}

\begin{proof}
  Let $b = b_1 \otimes \cdots \otimes b_p \in \hw_{I_0}^{\le \ell_0}(B)$, and 
  define $b_j^{(i)} \in B_j$ for $1 \le i \le j \le p$ as Proposition \ref{Prop:energy_function} (ii).
  Note that, since the combinatorial $R$-matrix and the action of $\tau$ commute,
  the first tensor factor of the image of $\tau(b_i \otimes \cdots \otimes b_j)$ under the isomorphism
  \[ B_i \otimes \cdots \otimes B_{j-1} \otimes B_j \stackrel{\sim}{\to} B_j \otimes B_{i} \otimes \cdots \otimes B_{j-1}
  \]
  is $\tau(b_j^{(i)})$.
  For each $1\le j \le p$, since $b \in \hw_{I_0}^{\le \ell_0}(B)$ implies $b_j^{(1)} \in \hw_{I_0}^{\le \ell_0}(B_j)$ 
  by  (\ref{eq:tensor}), we have that
  \[ \gs_{B_0, B_j} \big(m(B_0) \otimes \tau(b_j^{(1)})\big) =  b_j^{(1)} \otimes b^j
  \]
  for some $b^j \in B_0$ by Lemma \ref{Lem:Lemma2}.
  Hence by Proposition \ref{Prop:energy_function} (ii), we have
  \begin{align} \label{eq:D}
    D\big(m(B_0)\otimes \tau(b)\big) = D\big(m(B_0)\big) &+ \sum_{1 \le j \le p} D(b_j^{(1)}) + \sum_{1 \le j \le p} 
    H\big(m(B_0) \otimes \tau(b_j^{(1)})\big)  \nonumber \\
                                   & + \sum_{1 \le j < k \le p} H\big(\tau(b_j) \otimes \tau(b_k^{(j+1)})\big).  
  \end{align}
  For each $1 \le j \le p$ we have by Lemma \ref{Lem:Lemma3} that 
  \[ H\big(m(B_0) \otimes \tau(b_j^{(1)}) \big) = -\big\langle \wt(b_j^{(1)}), \varpi_{\tau^{-1}(0)}^\vee \big\rangle + C_j
  \]
  with some constant $C_j$ independent of $b_j^{(1)}$, and for each $1 \le j < k \le p$ we have by Lemma \ref{Lem:Lemma1} that
  \[ H\big(\tau(b_j) \otimes \tau(b_k^{(j+1)})\big) = H\big(b_j \otimes b_k^{(j+1)}\big) - \big\langle \wt(b_k^{(j+1)}) - 
     \wt(b_k^{(j)}), \varpi_{\tau^{-1}(0)}^\vee \big\rangle.
  \]
  Hence, it follows with some global constant $C$ that
  \begin{align*} (\ref{eq:D}) &= \sum_{1 \le j \le p} D(b_j^{(1)}) + \sum_{1\le j < k \le p} 
                 H\big(b_j \otimes b_k^{(j+1)}\big) - 
                 \sum_{1 \le j \le p} \big\langle \wt(b_j^{(1)}), \varpi_{\tau^{-1}(0)}^\vee \big\rangle \\
                 &\hspace{73pt}- \sum_{1 \le j < k \le p}
                 \big\langle \wt(b_k^{(j+1)}) - \wt(b_k^{(j)}), \varpi_{\tau^{-1}(0)}^\vee \big\rangle + C \\
                 &= D(b) - \sum_{1 \le k \le p} \big\langle \wt(b_k), \varpi_{\tau^{-1}(0)}^\vee \big\rangle + C \\
                 &= D(b) - \big\langle \wt(b), \varpi_{\tau^{-1}(0)}^\vee \big\rangle + C.
  \end{align*} 
  The assertion is proved.               
\end{proof}

Now we give the proof of our main theorem:\\

\noindent \textit{Proof of Theorem \ref{Thm:Main_Theorem}.\ }
  Let $\Psi_B$ be the isomorphism constructed in Proposition \ref{Prop:Main_Prop}, 
  and set $\ol{D}(b) = D(b) + \langle \wt\, \Psi_B (u_{\ell_p\gL_0} \otimes b),d \rangle$ for $b \in B$.
  It remains to verify that $\ol{D}(b) = C_B$ for every $b \in B$.
  First we show the following claim.\\
  \\
  \textit{Claim}. \ \ If there is some constant $C$ such that $\ol{D}(b) = C$ for all $b \in \hw_{I_0}^{\le \ell_p}(B)$,
  then $C = C_B$, and $\ol{D}(b) = C_B$ holds for all $b \in B$. \\[-0.2cm]
  
  By the second assertion of Proposition \ref{Prop:Main_Prop} and (\ref{eq:0}), we have
  \[ \ol{D}(u(B)) = \sum_{j=1}^p \big\langle t_{\mu_p + \cdots + \mu_j}(\ell^j\gL_0), d\big\rangle =C_B. 
  \]
  Hence it suffices to show that $\ol{D}(b) = C$ for all $b \in B$.
  Since $u_{\ell_p \gL_0} \otimes B$ is isomorphic to a disjoint union of Demazure crystals (see Remark \ref{Rem} (ii)),
  there exists a sequence $i_1,\ldots,i_k$ of elements of $I$ such that
  $e_{i_k} \cdots e_{i_1}(u_{\ell_p\gL_0} \otimes b)$ is a nonzero $U_q(\fg)$-highest weight element.
  We show $\ol{D}(b) = C$ by the induction on $k$.
  If $k =0$, then $u_{\ell_p\gL_0} \otimes b$ itself is $U_q(\fg)$-highest weight, 
  which is equivalent to $b \in \hw_{I_0}^{\le \ell_p}(B)$.
  Hence there is nothing to prove.
  Assume $k > 0$, and set $b' = e_{i_1}(b)$.
  Note that if $i_1 = 0$, then $\gee_0(b) > \ell_p$ holds since $e_0(u_{\ell_p\gL_0} \otimes b) = u_{\ell_p\gL_0} \otimes e_0b$.
  Hence we have $D(b') = D(b) - \gd_{0i_1}$ by Lemma \ref{Lem:elementary_for_energy}.
  On the other hand, it follows that
  \begin{align*}
    \wt\,\Psi_B(u_{\ell_p \gL_0} \otimes b') &= \wt\circ e_{i_1}\big(\Psi_B(u_{\ell_p \gL_0} \otimes b)\big)\\
                                             &= \wt\,\Psi_B(u_{\ell_p \gL_0} \otimes b) + \ga_{i_1}.
  \end{align*}
  Therefore we have
  \begin{align*}
     \ol{D}(b')= \big(D(b) -\gd_{0i_1}\big) + \Big(\big\langle \wt\, \Psi_B(u_{\ell_p \gL_0} \otimes b), d\big\rangle + \gd_{0i_1}\Big) = \ol{D}(b),
  \end{align*}
  and $\ol{D}(b) = C$ follows by the induction hypothesis. The claim is proved.\\
  
  In particular, since the set $\hw_{I_0}^{\le \ell_1}(B_1)$ contains only a single element $m(B_1)$,
  the theorem for $p=1$ follows from the claim.
  Assume $p>1$. We show the theorem by the induction on $p$. 
  Put $B^{p-1} = B_{p-1} \otimes \cdots \otimes B_1$ and $\tau= \tau^{r_p}$.
  Let $b \in \hw_{I_0}^{\le \ell_p}(B)$ be an arbitrary element, and write $b = b_p \otimes b^{p-1} \in B_p \otimes B^{p-1}$.
  Since $\lev(B_p)\ge \lev(B^{p-1})$, $b \in \hw_{I_0}^{\le \ell_p}(B)$ implies by Lemma \ref{Lem:level} (ii) that 
  \[ b_p \in \hw_{I_0}^{\le \ell_p}(B_p) \ \ \  \text{and} \ \ \  \gee(b^{p-1}) \in \gph\big(b_p\big) - P_{\cl}^+.
  \]
  Since $\lev(B_p) = \ell_p$, these are equivalent to 
  \[ b_p = m(B_p) \ \ \ \text{and} \ \ \ \gee(b^{p-1}) \in \ell_p \gL_{\tau(0)} - P_\cl^+.
  \]
  Hence if we put $b' = \tau^{-1}(b^{p-1})$, we have 
  \[ b = m(B_p) \otimes \tau(b')  \ \ \text{with} \ \ b' \in \hw_{I_0}^{\le \ell_p}(B^{p-1}).
  \]
  We see from the diagram in Remark \ref{Rem} (i) that
  \begin{align*}
     \wt\, \Psi_B\big(u_{\ell_p\gL_0} \otimes m(B_p) \otimes \tau(b')\big) &= \wt \circ \wti{\tau}\big(u_{\ell^p\gL_0}\otimes 
           \Psi_{B^{p-1}}(u_{\ell_{p-1}\gL_0} \otimes b')\big)\\
           &= \ell^p\gL_{\tau(0)} + \tau\circ \wt\, \Psi_{B^{p-1}}(u_{\ell_{p-1}\gL_0} \otimes b').
  \end{align*}
  Put $t = \tau^{-1}(0)$ and $C= C_{B^{p-1}} - \ell_{p-1}\langle \gL_0, d\rangle$. 
  It follows that
  \begin{align*}
     \tau\circ \wt\, \Psi_{B^{p-1}}&(u_{\ell_{p-1}\gL_0} \otimes b') = \tau\Big(\aff\circ\wt(b') + \ell_{p-1}\gL_0 + \big(-D(b') + C\big)\gd\Big) \\
    &= \aff \circ \wt\big(\tau(b')\big) + \ell_{p-1} \gL_{\tau(0)} + \Big( -D(b') + ( \wt(b'), \varpi_t) 
       + C\Big)\gd \\
    &= \aff \circ \wt\big(\tau(b')\big) + \ell_{p-1} \gL_{\tau(0)} + \Big(-D(b') + \langle \wt(b'), \varpi_t^{\vee} \rangle + C\Big)\gd,
  \end{align*}
  where the first equality follows since $\Psi_{B^{p-1}}$ is a $U_q'(\fg)$-crystal isomorphism and 
  \[ \langle \wt \,\Psi_{B^{p-1}}(u_{\ell_{p-1}\gL_0} \otimes 
  b'), d \rangle = -D(b') + C_{B^{p-1}}
  \]
  by the induction hypothesis, the second one follows since 
  \[ \tau = (\tau^t)^{-1} = (t_{\varpi_t}w_t)^{-1} = w_t^{-1}t_{-\varpi_t}
  \]
  by Lemma \ref{Prop:special_auto} (iii), and the third one follows from (\ref{eq:final}).
  Hence we have
  \[ \langle \wt\, \Psi_B\big(u_{\ell_p\gL_0} \otimes m(B_p) \otimes \tau(b')\big), d \rangle = \ell_p \langle \gL_{\tau(0)}, d\rangle 
     -D(b') + \langle \wt(b'), \varpi_t^{\vee} \rangle + C.
   \]
  On the other hand, we have from Lemma \ref{Lem:important} that
  \[ D\big(m(B_p) \otimes \tau(b')\big) = D\big(b'\big) - \langle \wt(b'), \varpi_t^\vee \rangle + C'
  \]
  with some global constant $C'$.
  Hence it is proved for all $ b \in \hw_{I_0}^{\le \ell_p}(B) = m(B_p) \otimes \Big\{ \tau(b') \Bigm| b' \in \hw_{I_0}^{\le \ell_p}(B^{p-1})\Big\}$ that
  \[ \ol{D}(b) = D(b) + \langle \wt\, \Psi_B (u_{\ell_p\gL_0} \otimes b),d \rangle = \ell_p \langle \gL_{\tau(0)}, d\rangle + C + C'.
  \]  
  Now the theorem follows from the claim. \qed

\def\cprime{$'$}


\end{document}